\RequirePackage{fix-cm}
\documentclass[smallextended]{svjour3}       
\smartqed  
\usepackage{graphicx}
\usepackage{amsmath}
\usepackage{amssymb}
\usepackage{latexsym}
\usepackage{mathrsfs}
\usepackage{colortbl}
\usepackage{amscd}
\usepackage{tikz-cd}
\usepackage{comment}

 \journalname{Japan Journal of Industrial and Applied Mathematics}

\def\div{\mathop{\mathrm{div}}}
\def\diam{\mathop{\mathrm{diam}}}

\def\>{\textgreater}
\def\<{\textless}
\def\Span{\mathop{\mathrm{span}}}

\spnewtheorem{thr}{Theorem}{\bf}{\it}
\spnewtheorem{defi}{Definition}{\bf}{\it}
\spnewtheorem{lem}{Lemma}{\bf}{\it}
\spnewtheorem{coro}{Corollary}{\bf}{\it}
\spnewtheorem{ass}{Assumption}{\bf}{\it}
\spnewtheorem*{pf*}{Proof}{\bf}{\rm}
\spnewtheorem*{rem*}{Remark:}{\it}{\it}
\spnewtheorem*{ex*}{Example:}{\it}{\it}
\spnewtheorem{rem}{Remark}{\bf}{\it}

\spnewtheorem*{lem1*}{Lemma 1}{\bf}{\rm}

\newcounter{sone}
\newcounter{stwo}
\newcounter{sthree}
\newcounter{sfour}
\newcounter{sfive}
\newcounter{ssix}
\newcounter{lone}
\newcounter{ltwo}
\newcounter{lthree}
\newcounter{lfour}
\newcounter{lfive}
\newcounter{lsix}
\makeatletter
    
    \@addtoreset{equation}{section}
  \makeatother
  
\begin{document}

\title{Crouzeix--Raviart and Raviart--Thomas finite-element error analysis on anisotropic meshes violating the maximum-angle condition
}

\titlerunning{Error analysis of the CR finite-element method}        

\author{Hiroki Ishizaka \and Kenta Kobayashi \and Takuya Tsuchiya 
}


\institute{Hiroki Ishizaka \at
              Graduate School of Science and Engineering, Ehime University, Matsuyama, Japan \\
              \email{h.ishizaka005@gmail.com}           
           \and
           Kenta Kobayashi\at
              Graduate School of Business Administration, Hitotsubashi University, Kunitachi, Japan \\
            \email{kenta.k@r.hit-u.ac.jp}
            \and
            Takuya Tsuchiya \at 
             Graduate School of Science and Engineering, Ehime University, Matsuyama, Japan \\
              \email{tsuchiya@math.sci.ehime-u.ac.jp}  
}

\date{Received: date / Accepted: date}

\maketitle

\begin{abstract}
We investigate the piecewise linear nonconforming Crouzeix--Raviar and the lowest order Raviart--Thomas finite-element methods for the Poisson problem on three-dimensional anisotropic meshes. We first give error estimates of the Crouzeix--Raviart and the Raviart--Thomas finite-element approximate problems. We next present the equivalence between the Raviart--Thomas finite-element method and the enriched Crouzeix--Raviart finite-element method. We emphasise that we do not impose either shape-regular or maximum-angle condition during mesh partitioning. Numerical results confirm the results that we obtained.




\end{abstract}

\section{Introduction}
\label{intro}
Let $\Omega \subset \mathbb{R}^d$, $d \in \{2,3\}$, be a bounded polyhedral domain. Furthermore, we assume that $\Omega$ is convex if necessary. We consider the Poisson problem as follows. Find $u: \Omega \to \mathbb{R}$ such that
\begin{align}
\displaystyle
- \varDelta u = f \quad \text{in $\Omega$}, \quad u = 0 \quad \text{on $\partial \Omega$}, \label{intro1}
\end{align}
where $f \in L^2(\Omega)$ is a given function. This paper gives error estimates for the first-order Crouzeix--Raviart (CR) finite-element approximation on anisotropic meshes in three dimensions. Anisotropic meshes have different mesh sizes in different directions. The shape regularity assumption on triangulations $\mathbb{T}_h$ is no longer valid on these meshes; see for example \cite{Ape99}. Furthermore, we do not impose the maximum-angle condition proposed in \cite{BabAzi76} during mesh partitioning. In many instances, the discussion also relates to two dimensions. We therefore discuss the problem here as uniformly valid in an arbitrary number of dimensions. 

CR finite error estimates for the non-homogeneous Dirichlet Poisson problem are known. Let $CR_{h0}^1$ be the CR finite-element space, to be defined in Section \ref{FEspaces}. Let $u \in H_0^1(\Omega)$ and $u_h^{CR} \in CR_{h0}^1$ be the exact and CR finite-element solutions, respectively. In \cite[Corollary 2.2]{Gud10}, adopting medius analysis, the estimate
\begin{align}
\displaystyle
| u - u_h^{CR} |_{H^1(\mathbb{T}_h)} 
\leq c_0 \left( \inf_{v_h \in CR_{h0}^1} | u - v_h |_{H^1(\mathbb{T}_h)} +  Osc_1(f) \right), \label{intro2}
\end{align}
is given, where $| \cdot |_{H^1(\mathbb{T}_h)}$ denotes the broken (piecewise) $H^1$-semi norm defined in Section \ref{regularmesh}, and $c_0$ a positive constant independent of $h$. Here, the oscillation $Osc_1(f)$ is expressed as
\begin{align*}
\displaystyle
Osc_1(f) := \left( \sum_{T \in \mathbb{T}_h} h_T^2 
\left [ \inf_{\bar{f} \in \mathcal{P}^0(T)} \| f - \bar{f} \|^2_{L^2(T)} \right] \right)^{1/2},
\end{align*}
where $\mathcal{P}^0(T)$ denotes the space of constant functions on $T$. Suppose that $u \in H^2(\Omega)$ and oscillation $Osc_1(f)$ vanishes. Let $I_h u \in CR_{h0}^1$ be the nodal interpolation of $u$ at the midpoints of the faces. Then, from the standard interpolation error estimate (see for example \cite[Corollary 1.109]{ErnGue04}), we have
\begin{align*}
\displaystyle
| u - u_h^{CR} |_{H^1(\mathbb{T}_h)} \leq c_0 | u - I_h u |_{H^1(\mathbb{T}_h)} 
\leq c_1 h |u|_{H^2(\Omega)}, 
\end{align*}
where $c_1$ represents a positive constant independent of $h$ and $u$ but depending on the parameter of the simplicial mesh; see for example \cite[Definition 1.107]{ErnGue04}. This parameter is bounded if the simplicial mesh sequence is shape regular. However, the situation is different without the shape-regular condition. The aim of the present paper is to deduce an analogous error estimate on anisotropic finite-element meshes. Note that very flat elements might be included in the mesh sequence. In many papers reporting on such investigations, the maximum-angle condition instead of the shape-regular condition is imposed. However, the maximum-angle condition is not necessarily needed to obtain error estimates. Recently, in the two-dimensional instance, the CR finite-element analysis of the non-homogeneous Dirichlet-Poisson problem has been investigated under a more relaxed mesh condition, \cite{KobTsu18b}. The present paper extends previous research to a three-dimensional setting.

However, it may not be easy to use the estimate \eqref{intro2} on anisotropic finite-element meshes. To overcome this difficulty, we use the interpolation error estimates obtained in \cite{IshKobTsu}. In that paper, the CR and Raviart--Thomas (RT) interpolation errors are bounded in terms of $h$ and the new parameter $H$, see Corollary \ref{modify=coro2}, \ref{modify=coro3}.

The CR finite-element space is not in $H_0^1(\Omega)$. Hence, an error between the exact solution and the CR finite-element approximation solution with a $H^1$-broken seminorm is divided into two parts (\cite{Bre15,ErnGue04}). One is an approximation error that measures how well the exact solution is approximated by the CR finite-element functions, the other is a nonconformity error term. For the former, the CR interpolation error estimates (Corollary \ref{modify=coro2}) are used. In the latter, the standard scaling argument is often used to obtain the error estimates. However, in this way, we are unable to derive the correct order on anisotropic meshes. To overcome this difficulty, we shall use the lowest-order RT interpolation error estimates on anisotropic meshes (Corollary \ref{modify=coro3}). By this technique, we consequently have the error estimates in the $H^1$-broken seminorm (Theorem \ref{CR=thr4}) and the $L^2$ norm (Theorem \ref{CR=thr5}) on anisotropic meshes. 

Furthermore, we present an error estimate for the first-order RT finite-element approximation of the Poisson problem \eqref{intro1} based on the dual mixed formulation (Theorem \ref{mix=thr7}). In the proof, we again use Corollary \ref{modify=coro3}. We again emphasise that we do not impose either the shape-regular or the maximum-angle condition during mesh partitioning.

We next present the equivalence of the enriched piecewise linear CR finite-element method introduced by \cite{HuMa15} and the first-order RT finite-element method. In two dimensions, the work \cite{ArnBre85} represents pioneering research. Marini \cite{Mar85} further found an expression relating RT and CR finite-element methods:
\begin{align}
\displaystyle
\bar{\sigma}_h^{RT}|_T &= \nabla \bar{u}_h^{CR} - \frac{f_T^0}{2} (x - x_T) \quad \text{on $T$}, \label{intro3}  
\end{align}
where $T$ denotes a mesh element, $x_i$ ($i=1,2,3$) the vertices of triangle $T$, $x_T$ the barycentre of $T$ such that $x_T := \frac{1}{3}(x_1 + x_2 + x_3)$, and $\bar{\sigma}_h^{RT}$ and $ \bar{u}_h^{CR}$ respectively denote the RT and CR finite-element solutions with a given external piecewise-constant function $f_T^0$. It was recently proved \cite{HuMa15} that the enriched piecewise-linear CR finite-element method is identical to the first-order RT finite-element method for both the Poisson and Stokes problems in any number of dimensions. In the present paper, we extend Marini's results to three dimensions (Lemma \ref{RTCR=lem12}). 

The remainder of the present paper is organised as follows. Section 2 introduces the weak form of the continuous problem \eqref{intro1}, the finite-element meshes, and finite-element spaces. Furthermore, we propose a parameter $H$. Section 3 introduces discrete settings of the CR finite-element method for \eqref{intro1} and proposes error estimates. Section 4 proves error estimates for the first-order RT finite-element method based on the dual mixed formulation of the Poisson problem. Section 5 gives the equivalence of the RT and CR finite-element problems. Finally, Section 6 presents numerical results obtained using the Lagrange P1 element and the first-order CR element.

\section{Preliminaries}
\subsection{Weak formulation}
The variational formulation for the Poisson problem \eqref{intro1} is then as follows. Find $u \in H_0^1(\Omega)$ such that
\begin{align}
\displaystyle
a_0(u,\varphi) = (f , \varphi) \quad \forall \varphi \in H_0^1(\Omega), \label{pre1} 
\end{align}
where $a_0: H^1(\Omega) \times H^1(\Omega) \to \mathbb{R}$ denotes a bilinear form defined by
\begin{align*}
\displaystyle
a_0(u,\varphi) := (\nabla u , \nabla \varphi).
\end{align*}
Here, we define  $H_0^1(\Omega)$ as the closure of $C_0^{\infty}(\Omega)$ in the semi-norm  $| \cdot |_{H^1(\Omega)}$. By the Lax--Milgram lemma, there exists a unique solution $u \in H_0^1(\Omega)$ for any $f \in L^2(\Omega)$ and it holds that
\begin{align*}
\displaystyle
| u |_{H^1(\Omega)} \leq C_P(\Omega) \| f \|,
\end{align*}
where $C_P(\Omega)$ is the Poincar$\rm{\acute{e}}$ constant depending on $\Omega$. 
Furthermore, if $\Omega$ is convex, then $u \in H^2(\Omega)$ and 
\begin{align}
\displaystyle
| u |_{H^2(\Omega)} \leq \| \varDelta u \|. \label{pre2}
\end{align}
The proof can be found in, for example, \cite[Theorem 3.1.1.2, Theorem 3.2.1.2]{Gri11}. 

\subsection{Meshes, Mesh faces, Averages and Jumps} \label{regularmesh}
Let $\mathbb{T}_h = \{ T \}$ be a simplicial mesh of $\overline{\Omega}$, made up of closed $d$-simplices, such as
\begin{align*}
\displaystyle
\overline{\Omega} = \bigcup_{T \in \mathbb{T}_h} T,
\end{align*}
with $h := \max_{T \in \mathbb{T}_h} h_{T}$, where $ h_{T} := \diam(T)$. We assume that each face of any $d$-simplex $T_1$ in $\mathbb{T}_h$ is either a subset of the boundary $\partial \Omega$ or a face of another $d$-simplex $T_2$ in  $\mathbb{T}_h$. That is, $\mathbb{T}_h$ is a simplicial mesh of $\overline{\Omega}$ without hanging nodes. 

\begin{defi} \label{pre=defi1}
For any $T \in \mathbb{T}_h$, we define the parameter $H_{T}$ as
\begin{align*}
\displaystyle
H_{T} := \frac{h_{T}^2}{|T|} \min_{1 \leq i \leq 3} |L_i|  \quad \text{if $d=2$},
\end{align*}
where $L_i$ $(i=1,2,3)$ denotes edges of the triangle $T$. Further, we define the parameter $H_{T}$ as
\begin{align*}
\displaystyle
H_{T} := \frac{h_{T}^2}{|T|} \min_{1 \leq i , j \leq 6, i \neq j} |L_i| |L_j| \quad \text{if $d=3$},
\end{align*}
where $L_i$ $(i=1,\ldots,6)$ denotes edges of the tetrahedra $T$. Here, $|T|$ denotes the measure of $T$. Furthermore, we set 
\begin{align*}
\displaystyle
H := H(h) := \max_{T \in \mathbb{T}_h} H_{T}.
\end{align*}
\end{defi}
We impose the following assumption.
 \begin{ass}  \label{pre=ass1}
 We assume that $\{ \mathbb{T}_h \}_{h \> 0}$ is a sequence of  triangulations of $\Omega$ such that
\begin{align*}
\displaystyle
\lim_{h \to 0} H(h) = 0.
\end{align*}
 \end{ass}
 
 \begin{rem}
 The parameter $H_T$ was introduced, and the interpolation errors are bounded (locally) in terms of $h_T$ and $H_T$ on anisotropic meshes without any geometric conditions in \cite{IshKobTsu}. In two-dimensional case, the parameter $H_T$ is equivalent to the circumradius of $T$.  Hence, the maximum-angle condition or the semiregular condition holds if and only if there exists a constant $\sigma_0$ such that $H_T / h_T \leq \sigma_0$. In three-dimensional case, it is conjectured that the maximum-angle condition holds if and only if the quantity $H_T / h_T$ is bounded.
 \end{rem}

We adopt the concepts of mesh faces, averages and jumps in the analysis of RT and CR finite element method. Let $\mathcal{F}_h^i$ be the set of interior faces and $\mathcal{F}_h^{\partial}$ the set of the faces on the boundary $\partial \Omega$. Let $\mathcal{F}_h := \mathcal{F}_h^i \cup \mathcal{F}_h^{\partial}$. For any $F \in \mathcal{F}_h$, we define the unit normal $n_F$ to $F$ as follows: (\roman{sone}) If  $F \in \mathcal{F}_h^i$ with $F = T_1 \cap T_2$, $T_1,T_2 \in \mathbb{T}_h$, let $n_1$ and $n_2$ be the outward unit normals of $T_1$ and $T_2$, respectively. Then, $n_F$ is either of $\{ n_1 , n_2\}$; (\roman{stwo}) If $F \in \mathcal{F}_h^{\partial}$, $n_F$ is the unit outward normal $n$ to $\partial \Omega$.

Let $k$ be a positive integer. We then define the broken (piecewise) Sobolev space as
\begin{align*}
\displaystyle
H^k(\mathbb{T}_h) &:= \left\{ \varphi \in L^2(\Omega); \ \varphi|_T \in H^k(T) \ \forall T \in \mathbb{T}_h  \right\}
\end{align*}
with the norm
\begin{align*}
\displaystyle
| \varphi |_{H^1(\mathbb{T}_h)} &:= \left( \sum_{T \in \mathbb{T}_h}\| \nabla \varphi \|^2_{L^2(T)^d} \right)^{1/2} \quad \varphi \in H^1(\mathbb{T}_h).
\end{align*}

Let $\varphi \in H^k(\mathbb{T}_h)$. Suppose that $F \in \mathcal{F}_h^i$ with $F = T_1 \cap T_2$, $T_1,T_2 \in \mathbb{T}_h$. Set $\varphi_1 := \varphi{|_{T_1}}$ and $\varphi_2 := \varphi{|_{T_2}}$. The jump and the average of $\varphi$ across $F$ is then defined as
\begin{align*}
\displaystyle
 [[ \varphi ]]_F := (\varphi_1 n_1 + \varphi_2 n_2) \cdot n_F, \quad \{ \! \{ \varphi \} \! \}_F := \frac{1}{2} (\varphi_1 + \varphi_2).
\end{align*}
For a boundary face $F \in \mathcal{F}_h^{\partial}$ with $F = \partial T \cap \partial \Omega$, $[[\varphi]]_F := \varphi|_T$ and $\{ \! \{ \varphi \} \! \}_F := \varphi |_T$. When $v$ is an $\mathbb{R}^d$-valued function, we use the notation
\begin{align*}
\displaystyle
[[ v \cdot n ]]_F := (v_1 - v_2) \cdot n_F, \quad \{ \! \{ v\} \! \}_F := \frac{1}{2} (v_1  + v_2 )
\end{align*}
for the jump of the normal component of $v$. For a boundary face $F \in \mathcal{F}_h^{\partial}$ with $F = \partial T \cap \partial \Omega$, $[[v \cdot n]]_F := v|_T \cdot n $ and $\{\! \{ v  \}\! \}_F := v|_T $. Whenever no confusion can arise, we simply write $[[\varphi]]$, $\{\! \{ \varphi \}\! \}$, $[[v \cdot n]]$ and $\{\! \{ v \} \! \}$, respectively.

 Suppose that $F \in \mathcal{F}_h^i$ with $F = T_1 \cap T_2$, $T_1,T_2 \in \mathbb{T}_h$. For $v \in H^1(\mathbb{T}_h)^d$ and $\varphi \in H^1(\mathbb{T}_h)$, it holds that
\begin{align*}
\displaystyle
[[ (v \varphi) \cdot n ]]_F = \{ \! \{ v \} \! \}_F \cdot n_F [[ \varphi ]]_F + [[ v \cdot n]]_F  \{ \! \{ \varphi \} \! \}_F.
\end{align*}

We here define a broken gradient operator as follows.
\begin{defi} \label{pre=defi2}
For $\varphi \in H^1(\mathbb{T}_h)$, the broken gradient $\nabla_h:H^1(\mathbb{T}_h) \to L^2(\Omega)^d$ is defined by
\begin{align*}
\displaystyle
(\nabla_h \varphi)|_T &:= \nabla (\varphi|_T) \quad \forall T \in \mathbb{T}_h.
\end{align*}
\end{defi}
\noindent
Note that $H^1(\Omega )\subset H^1(\mathbb{T}_h)$ and the broken gradient coincides with the distributional gradient in $H^1(\Omega)$. 

\subsection{Finite Element Spaces and Interpolations Error Estimates} \label{FEspaces}
This section introduce the piecewise-constant, CR and RT finite element spaces.

Let $T \in \mathbb{T}_h$. For any $k \in \mathbb{N}_0$, let $\mathcal{P}^k(T)$ be the space of polynomials with degree at most $k$ in $T$. 


\begin{thr}[Poincar\'e inequality] \label{modify=thr1}
Let $D \subset \mathbb{R}^d$  be a convex domain with diameter $\diam(D)$. It then holds that, for $\varphi \in H^1(D)$ with $\int_D \varphi dx = 0$,  
\begin{align}
\displaystyle
\| \varphi \|_{L^2(D)} \leq \frac{\diam (D)}{\pi} |\varphi|_{H^1(D)}.\label{modify3}
\end{align}
\end{thr}

\begin{pf*}
The proof is found in \cite[Theorem 3.2]{Mar03}, also see \cite{PayWei60}.
\qed
\end{pf*}

\subsubsection{Piecewise-constant finite element space}
We define the standard piecewise constant space as
\begin{align*}
\displaystyle
M_h^0 := \left\{ q_h \in L^2(\Omega); \ q_h|_T \in \mathcal{P}^0(T) \ \forall T \in \mathbb{T}_h \right\}.
\end{align*}
The local interpolation $\Pi_{T}^0$ from $L^2(T)$ into the space $\mathcal{P}^0(T)$ is defined by
\begin{align*}
\displaystyle
\int_{T} (\Pi_{T}^0 q - q) dx = 0 \quad \forall q \in L^2(T).
\end{align*}
Note that $\Pi_{T}^0 q$ is the constant function equal to $\frac{1}{|T|} \int_T q dx$. We also define the global interpolation $\Pi_h^0$ to the space $M_{h}^0$ by
\begin{align*}
\displaystyle
(\Pi_h^0 q)|_T = \Pi_T^0 (q|_T) \quad \forall T \in \mathbb{T}_h, \quad \forall q \in L^2(\Omega). 
\end{align*}

The Poincar\'e inequality \eqref{modify3} directly yields the following error estimate of the local $L^2$-projection $\Pi^0_T$.

\begin{thr}  \label{modify=thr2}
We have the error estimate of the local $L^2$-projection such that
\begin{align}
\displaystyle
\| \Pi^0_T q - q \|_{L^2(T)} &\leq \frac{h_T}{\pi} |q|_{H^1(T)} \quad \forall T \in \mathbb{T}_h, \quad \forall q \in H^1(T). \label{modify6}
\end{align}
\end{thr}

\begin{pf*}
For any $q \in H^1(T)$, we set $w := \Pi^0_T q - q$. It then holds that
\begin{align*}
\displaystyle
\int_T w dx = \int_T ( \Pi^0_T q - q ) dx = \frac{1}{|T|} \int_T q dx |T| - \int_T q dx = 0.
\end{align*}
Therefore, using the Poincar\'e inequality \eqref{modify3}, we conclude \eqref{modify6}.
\qed	
\end{pf*}

The global error estimate of the $L^2$-projection is obtained as follows.
\begin{coro}  \label{modify=coro1}
Let $\{ \mathbb{T}_h \}$ be a family of conformal meshes satisfying Assumption \ref{pre=ass1}. It then holds that
\begin{align}
\displaystyle
\| \Pi^0_h q - q \|_{L^p(\Omega)} &\leq \frac{h}{\pi} |q|_{H^{1}(\Omega)} \quad \forall q \in H^{1}(\Omega). \label{modify7}
\end{align}
\end{coro}

\subsubsection{CR finite element space}
We define the following CR finite element space as
\begin{align*}
\displaystyle
CR_{h0}^1 &:= \left \{ \varphi_h \in L^2(\Omega); \ \varphi_h|_T \in \mathcal{P}^1(T) \ \forall T \in \mathbb{T}_h, \ \int_F [[ \varphi_h ]]_F ds = 0 \ \forall F \in \mathcal{F}_h \right \}.
\end{align*}
Using the barycentric coordinates $ \lambda_i: \mathbb{R}^d \to \mathbb{R}$, $i=1,\ldots,d+1$, we define the local basis functions as
\begin{align*}
\displaystyle
\theta_i(x) := d \left( \frac{1}{d} - \lambda_i(x) \right), \quad 1 \leq i \leq d+1. 
\end{align*}
For $i=1,\ldots,d+1$, let $F_i$ be the face of $T$ and $x_{F_i}$ the barycentre of the face $F_i$. We then define the local CR interpolation operator as
\begin{align*}
\displaystyle
I_{T}^{CR}: H^{1}(T) \ni \varphi  \mapsto I_{T}^{CR} \varphi := \sum_{i=1}^{d+1} \left( \frac{1}{|F_i|} \int_{F_i} \varphi ds \right)  \theta_i \in \mathcal{P}^1. 
\end{align*}
Furthermore, it holds that
\begin{align*}
\displaystyle
 \frac{1}{|F_i|} \int_{F_i} \left( I_{T}^{CR} \varphi -  \varphi \right) ds = 0, \quad i=1,\ldots,d+1, \quad \forall \varphi \in H^1(T).
\end{align*}
We define the global CR interpolation $I_h^{CR}:H_0^{1}(\Omega) \to CR_{h0}^1$ by
\begin{align*}
\displaystyle
(I_h^{CR} \varphi )|_T = I_T^{CR} (\varphi|_T) \quad \forall T \in \mathbb{T}_h, \quad \forall \varphi \in H_0^{1}(\Omega).
\end{align*}

We give the local CR interpolation error estimate.

\begin{thr}  \label{modify=thr3}
We have the following estimates such that for $m \in \{ 0 , 1 \}$,
\begin{align}
\displaystyle
 |\varphi - {I}^{CR}_{T} \varphi|_{H^{m}(T)}
&\leq C_I^{CR,m} \left( \frac{H_{T}}{h_{T}} \right)^m h_{T}^{2-m} | \varphi |_{H^{2}(T)} \ \forall T \in \mathbb{T}_h, \ \forall \varphi \in H^{2}(T). \label{modify10}
\end{align}
Here,  $C_I^{CR,m}$ is a positive constant independent of $h_T$ and $H_T$.
\end{thr}

\begin{pf*}
The proof is found in \cite[Theorem 2]{IshKobTsu}.
\qed
\end{pf*}

The global CR interpolation error estimates are obtained as follows.

\begin{coro}  \label{modify=coro2}
Let $\{ \mathbb{T}_h \}$ be a family of conformal meshes satisfying Assumption \ref{pre=ass1}. Then, there exists constants $C_G^{CR,0},  C_G^{CR,1}\> 0$, independent of $H$ and $h$, such that
\begin{align}
\displaystyle
 \|\varphi - {I}^{CR}_{h} \varphi \|_{L^{2}(\Omega)}
&\leq C_G^{CR,0} h^2 | \varphi |_{H^{2}(\Omega)} \quad  \forall \varphi \in H_0^1(\Omega) \cap H^{2}(\Omega),\label{modify11} \\
 |\varphi - {I}^{CR}_{h} \varphi|_{H^{1}(\mathbb{T}_h)}
&\leq C_G^{CR,1} H | \varphi |_{H^{2}(\Omega)} \quad  \forall \varphi \in H_0^1(\Omega) \cap H^{2}(\Omega). \label{modify12}
\end{align}
\end{coro}

The inequality \eqref{modify10} with $m=1$ can be improved by replacing $H_{T}$ with $h_{T}$. To this end, we use the Poincar\'e inequality \eqref{modify3}. 

\begin{thr}  \label{modify=thr4}
It then holds that
\begin{align}
\displaystyle
| I_{T}^{CR} \varphi - \varphi |_{H^1(T)}  &\leq \frac{h_{T}}{\pi} |\varphi|_{H^2(T)} \quad \forall T \in \mathbb{T}_h, \quad \forall \varphi \in H^2(T). \label{modify13}
\end{align}
\end{thr}

\begin{pf*}
Let $F_i$, $i = 1,\ldots,d+1$ be the faces of the element $T$. 
We set $\psi := I_{T}^{CR} \varphi - \varphi \in H^2(T)$. From Green's formula and the property of the CR interpolation, we have
\begin{align*}
\displaystyle
\int_{T} \frac{\partial \psi}{\partial x_j} dx = \int_{\partial T} \psi n_{T}^{(j)} ds = \sum_{i=1}^{d+1} n_{T}^{(j)} \int_{F_i} \psi ds = 0,
\end{align*}
where $ n_{T}^{(j)}$ denotes the $j$th component of the outer unit normal vector $n_{T}$. From  the Poincar\'e inequality \eqref{modify3}, we have
\begin{align*}
\displaystyle
| I_{T}^{CR} \varphi - \varphi |_{H^1(T)}^2
&= \sum_{j=1}^d \left \|  \frac{\partial}{\partial x_j} ( I_{T}^{CR} \varphi - \varphi)  \right\|^2_{L^2(T)} \\
&\leq \left( \frac{h_{T}}{\pi} \right)^2 \sum_{j=1}^d \left | \frac{\partial}{\partial x_j} ( I_{T}^{CR} \varphi - \varphi) \right |^2_{H^1(T)} \\
&=  \left( \frac{h_{T}}{\pi} \right)^2 \sum_{j,k=1}^d  \left \|  \frac{\partial^2}{\partial x_j \partial x_k} ( I_{T}^{CR} \varphi - \varphi) \right \|^2_{L^2(T)} \\
&= \left( \frac{h_{T}}{\pi} \right)^2 |\varphi|^2_{H^2(T)},
\end{align*}
which conclude \eqref{modify13}.
\qed
\end{pf*}

\begin{rem} \label{modify=rem1}
For $i=1,\ldots,d+1$, let $x_{F_i}$ the barycentre of face $F_i$. If we choose the domain of the local CR interpolation operator as $W^{\ell,p}(T) \subset \mathcal{C}^0(T)$ with $1 \leq p \< \infty$ and $d \< \ell p$, it is possible to define
\begin{align*}
\displaystyle
I_{T}^{CR,S}:W^{\ell,p}(T) \ni \varphi  \mapsto I_{T}^{CR,S} \varphi := \sum_{i=1}^{d+1} \varphi (x_{F_i}) \theta_i \in \mathcal{P}^1.
\end{align*}
However, the estimate \cite[Theorem 2]{IshKobTsu}
\begin{align*}
\displaystyle
| I_{T}^{CR,S} \varphi - \varphi |_{H^{1}(T)}  &\leq C_I^{CR,1} H_{T} |\varphi|_{H^{2}(T)}  \quad \forall \varphi \in H^{2}(T)
\end{align*}
can not be improved by replacing $H_{T}$ with $h_{T}$. 

%
%

As a counter example, let us consider $T$ with vertices $x_1 := (0,0,0)^T$, $x_2 := (h,0,0)^T$, $x_3 := (\frac{h}{2}, h^{\gamma},0)^T$ and  $x_4 := (\frac{h}{2}, 0, \frac{h}{2})^T$, where $h := \frac{1}{N}$, $N \in \mathbb{N}$ and $\gamma \in \mathbb{R}$, $1 \< \gamma \leq 2$. Let $\varphi$ be a function such that
\begin{align*}
\displaystyle
\varphi(x,y,z) := x^2 + y^2 + z ^2.
\end{align*}
If an exact solution $\varphi$ is known, the error $e_h := \varphi - \varphi_h$ and $e_{h/2} := \varphi - \varphi_{h/2}$ are computed numerically for two mesh sizes $h$ and $h/2$, where $\varphi_h := I_{T}^{CR,S} \varphi$. The convergence indicator $r$ is defined by
\begin{align*}
\displaystyle
r = \frac{1}{\log(2)} \log \left( \frac{\| e_h \|_X}{\| e_{h/2} \|_X} \right).
\end{align*}
The parameter $H_{T}$ is then $H_{T} = \mathcal{O}( h^{2 - \gamma} )$. We compute the convergence order with respect to the $H_0^1$ norm defined by
\begin{align*}
\displaystyle
&Err_h^{CR,S}(H^1) := \frac{| \varphi - I_{T}^{CR,S} \varphi |_{H^1(T)}}{| \varphi |_{H^2(T)}}, 
\end{align*}
for the case: $\gamma = 1.5$ (Table \ref{modify=table1}).

\begin{table}[htb]
\caption{Error of the local CR interpolation operator ($\gamma = 1.5$)}
\centering
\begin{tabular}{l | l | l | l | l } \hline
$N$ &  $h$ & $H_{T}$ & $Err_h^{CR,S}(H^1)$ & $r$  \\ \hline \hline
128 & 7.8125e-03 & 3.8081e-01  & 2.8183e-03  &    \\
256 & 3.9062e-03  & 2.6723e-01  &  1.7641e-03 & 0.68  \\
512 & 1.9531e-03  & 1.8823e-01  & 1.1587e-03   &  0.61  \\
1024 & 9.7656e-04  & 1.3284e-01  &  7.8625e-04  & 0.60  \\
2048 & 4.8828e-04 & 9.3842e-02  &  5.4390e-04 &  0.53 \\
4096 & 2.4414e-04  &  6.6324e-02 &  3.8026e-04  & 0.52 \\
\hline
\end{tabular}
\label{modify=table1}
\end{table}
\end{rem}

\subsubsection{RT finite element space}
The lowest order RT finite element space is defined by
\begin{align*}
\displaystyle
RT^0(T) := \{ v ; \ v(x) = p + x q, \ p \in \mathcal{P}^0(T)^d, \ q \in \mathcal{P}^0(T), \ x \in \mathbb{R}^d \}.
\end{align*}
The functionals are defined by, for any $v \in RT^0(T)$,
 \begin{align*}
\displaystyle
{\chi}_{i} (v) := \frac{1}{|F_i|} \int_{F_i} v \cdot n_{i} ds, \quad F_i \subset \partial T, \quad 1 \leq i \leq d+1, 
\end{align*}
where $n_i$ denotes the outer unit normal vector of $T$ along $F_i$. We set $\sum := \{ {\chi}_{i} \}_{i=1}^{d+1}$. Note that $\dim RT^0(T) = d+1$. The triple $\{ T , RT^0(T) , \Sigma \}$ is then a finite element. We define the RT finite element space by
\begin{align*}
\displaystyle
RT^0_{h} &:= \{ v_h \in L^2(\Omega)^d;  \ v_h |_T \in RT^0(T), \ \forall T \in \mathbb{T}_h, \ [[ v_h \cdot n ]]_F = 0, \ \forall F \in \mathcal{F}_h^i \}.
\end{align*}
Note that $RT^0_{h} \subset H(\div;\Omega)  := \left\{ v \in L^2(\Omega)^d ; \ \div v \in L^2(\Omega) \right\}$.

We next define the local RT interpolation as
\begin{align*}
\displaystyle
I_T^{RT}: H^1(T)^d \to RT^0(T), 
\end{align*}
using
\begin{align*}
\displaystyle
 \int_{F_i} (v - I_T^{RT} v) \cdot n_i ds = 0, \quad F_i \subset \partial T, \ i \in \{ 1 , \ldots, d+1\} \quad \forall v \in H^1(T)^d. 
\end{align*}
Further, we define the global RT interpolation $I_h^{RT} : H^1(\Omega)^d \to RT^0_{h}$ by
\begin{align*}
\displaystyle
(I_h^{RT} v )|_T = I_T^{RT} (v|_T) \quad \forall T \in \mathbb{T}_h, \quad \forall v \in H^1(\Omega)^d.
\end{align*}

The local RT interpolation error estimate is as follows.

\begin{thr}  \label{modify=thr5}
We have the following estimates such that
\begin{align}
\displaystyle
\| I_{T}^{RT} v - v \|_{L^2(T)^d}
&\leq  C_I^{RT} H_{T} |v|_{H^{1}(T)^d} \quad \forall T \in \mathbb{T}_h, \quad  \forall v \in H^{1}(T)^d. \label{modify14}
\end{align}
Here,  $C_I^{RT}$ is a positive constant independent of $H_T$.
\end{thr}

\begin{pf*}
The proof is found in \cite[Theorem 3]{IshKobTsu}.
\qed
\end{pf*}

The global RT interpolation error estimates are obtained as follows.

\begin{coro}  \label{modify=coro3}
Let $\{ \mathbb{T}_h \}$ be a family of conformal meshes satisfying Assumption \ref{pre=ass1}. Then, there exists a constant $C_G^{RT} \> 0$, independent of $H$, such that
\begin{align}
\displaystyle
\| I_{h}^{RT} v - v \|_{L^2(\Omega)^d}
&\leq  C_G^{RT} H |v|_{H^{1}(\Omega)^d} \quad  \forall v \in H^{1}(\Omega)^d. \label{modify15}
\end{align}
\end{coro}

Between the RT interpolation $I_h^{RT}$ and the $L^2$-projection $\Pi_h^0$, the following relation holds:
\begin{lem} \label{pre=lem1}
For any $v \in H^1(\Omega)^d$, it holds that
\begin{align*}
\displaystyle
\div (I_h^{RT} v) = \Pi_h^0 (\div v).
\end{align*}
That is to say, the diagram
\[
\begin{CD}
     H^1(\Omega)^d @>{\div}>> L^2(\Omega) \\
  @V{I_h^{RT}}VV    @VV{\Pi_h^0}V \\
     RT^0_{h}   @>{\div}>>  M_{h}^0
\end{CD}
 \]
commutes. 
\end{lem}

\begin{pf*}
The proof of this lemma is found in \cite{Bra07}.
\qed	
\end{pf*}

The following relation plays an important role in the CR finite element analysis on anisotropic meshes.
\begin{lem} \label{pre=lem2}
It holds that
\begin{align}
\displaystyle
(v_h , \nabla_h \psi_h) + (\div v_h , \psi_h) = 0 \quad \forall v_h \in RT_h^0, \quad \forall \psi_h \in H_0^1(\Omega) +  CR_{h0}^1.  \label{pre21}
\end{align}
\end{lem}

\begin{pf*}
For any $v_h \in RT_h^0$ and $\psi_h \in H_0^1(\Omega) +  CR_{h0}^1$, using Green formula and the fact $v_h \cdot n_F \in \mathcal{P}^0(F)$ for any $F \in \mathcal{F}_h$, we can derive 
\begin{align*}
\displaystyle
(v_h , \nabla_h \psi_h) + (\div v_h , \psi_h)
&= \sum_{T \in \mathbb{T}_h} \int_{\partial T} (v_h \cdot n_T) \psi_h ds  \notag \\
&= \sum_{F \in \mathcal{F}_h} \int_{F} [[ (v_h  \psi_h ) \cdot n_F ]] ds  \notag \\
&= \sum_{F \in \mathcal{F}_h} \int_{F} \left(  [[ v_h \cdot n_F]] \{ \! \{ \psi_h \} \!\} + \{ \! \{ v_h  \} \!\} \cdot n_F [[ \psi_h ]] \right) ds \notag \\
&=  0.
\end{align*}
\qed	
\end{pf*}

\subsection{Discrete Poincar\'e Inequality on Anisotropic Meshes}
We propose the discrete Poincar\'e inequality on anisotropic meshes. 
\begin{lem}[Discrete Poincar\'e inequality on anisotropic meshes] \label{pre=lem3}
Assume that $\Omega$ is convex.  If $H \leq 1$,  there exists $C(\Omega)$, independent of $h$, $H$, and the geometry of meshes, such that
\begin{align}
\displaystyle
 \| \varphi_h \| \leq C(\Omega) |\varphi_h|_{H^1(\mathbb{T}_h)} \quad \forall \varphi_h \in CR_{h0}^1. \label{pre22}
\end{align}
\end{lem}

\begin{pf*}
Let $\varphi_h \in CR_{h0}^1$. We consider the dual problem. Find $z \in H^{2}(\Omega) \cap H_0^1(\Omega)$ such that
\begin{align*}
\displaystyle
- \varDelta z = \frac{\varphi_h}{\| \varphi_h \|} \quad \text{in $\Omega$}, \quad z = 0 \quad \text{on $\partial \Omega$}.
\end{align*}
We then have a priori estimates:
\begin{align*}
\displaystyle
|z|_{H^1(\Omega)} \leq C_P, \quad |z|_{H^2(\Omega)} \leq 1,
\end{align*}
where $C_P$ is the Poincar\'e constant. We use the duality argument to show the target inequality. That is to say, we have
\begin{align*}
\displaystyle
\| \varphi_h \| &= \frac{1}{\| \varphi_h \|} (\varphi_h , \varphi_h) = (- \varDelta z , \varphi_h) = (- \div \nabla z , \varphi_h) \\
&= (- \div \nabla z , \varphi_h - \Pi_h^0 \varphi_h) - (\nabla z - I_h^{RT}(\nabla z) , \nabla_h \varphi_h) + (\nabla z , \nabla_h \varphi_h) \\
&\leq \| \varDelta z\| \|  \varphi_h - \Pi_h^0 \varphi_h \| + \| \nabla z - I_h^{RT}(\nabla z) \| |\varphi_h|_{H^1(\mathbb{T}_h)} + |z|_{H^1(\Omega)} |\varphi_h|_{H^1(\mathbb{T}_h)} \\
&\leq c \left( h+ H |\nabla z|_{H^1(\Omega)}+ C_P \right) |\varphi_h|_{H^1(\mathbb{T}_h)},
\end{align*}
which leads to
\begin{align*}
\displaystyle
\| \varphi_h \| 
\leq c (2 + C_p) |\varphi_h|_{H^1(\mathbb{T}_h)} \quad \text{if $H \leq 1$}.
\end{align*}
We here used
\begin{align*}
\displaystyle
- \int_{\Omega} \div (\nabla z)  \varphi_h dx
&=  \int_{\Omega} ( \Pi_h^0  \div (\nabla z) - \div (\nabla z) ) \varphi_h dx - \int_{\Omega} ( \Pi_h^0  \div (\nabla z) ) \varphi_h dx\\
&=  \int_{\Omega} ( \Pi_h^0  \div (\nabla z) - \div (\nabla z) ) ( \varphi_h - \Pi_h^0 \varphi_h) dx \\
&\quad  - \int_{\Omega} ( \div I_h^{RT}  (\nabla z) ) \varphi_h dx\\
&= - \int_{\Omega} \div (\nabla z) \left( \varphi_h - \Pi_h^0 \varphi_h \right) dx \\
&\quad - \int_{\Omega} (\nabla z - I_h^{RT} (\nabla z) ) \cdot \nabla_h \varphi_h dx + \int_{\Omega} \nabla z \cdot \nabla_h \varphi_h dx,
\end{align*}
where
\begin{align*}
\displaystyle
 \int_{\Omega} ( \div I_h^{RT}  (\nabla z) ) \varphi_h dx
 &= \sum_{T \in \mathbb{T}_h} \int_{\partial T} n_T \cdot  I_h^{RT}  (\nabla z) \varphi_h ds - \int_{\Omega}  I_h^{RT}  (\nabla z) \cdot \nabla_h \varphi_h dx \\
 &=  \int_{\Omega} ( \nabla z -  I_h^{RT}  (\nabla z) ) \cdot \nabla_h \varphi_h dx  - \int_{\Omega}  \nabla z \cdot \nabla_h \varphi_h dx.
\end{align*}
\qed	
\end{pf*}

\section{CR Finite Element Approximation}

\subsection{Finite Element Approximation}
The CR finite element problem is to find $u_h^{CR}  \in CR_{h0}^1$ such that
\begin{align}
\displaystyle
a_{0h}(u_h^{CR} , \varphi_h) = (f , \varphi_h) \quad \forall \varphi_h \in CR_{h0}^1, \label{CR1}
\end{align}
where $a_{0h}: ( CR_{h0}^1 + H_0^1(\Omega) ) \times ( CR_{h0}^1 + H_0^1(\Omega) ) \to \mathbb{R}$ is defined by
\begin{align*}
\displaystyle
a_{0h}(\psi_h,\varphi_h) := \sum_{T \in \mathbb{T}_h} \int_T \nabla \psi_h \cdot \nabla \varphi_h dx = (\nabla_h \psi_h , \nabla_h \varphi_h).
\end{align*}
This problem is nonconforming because $CR_{h0}^1 \not\subset H_0^1(\Omega)$. 

For the CR approximate solution $u_h^{CR}  \in CR_{h0}^1$ of \eqref{CR1}, we have the a priori estimate, using \eqref{pre22},
\begin{align*}
\displaystyle
| u_h^{CR} |_{H^1(\mathbb{T}_h)}^2
\leq \| f \| \| u_h^{CR} \|
\leq C(\Omega) \| f \| | u_h^{CR} |_{H^1(\mathbb{T}_h)}.
\end{align*}
By the Lax--Milgram lemma, there exists a unique solution $u_h^{CR} \in CR_{h0}^1$ for any $f \in L^2(\Omega)$.

\subsection{Classical Error Analysis}
The starting point for error analysis is the Second Strang Lemma, e.g. see \cite[Lemma 2.25]{ErnGue04},
\begin{align}
\displaystyle
| u - u_h^{CR} |_{H^1(\mathbb{T}_h)} \leq 2 \inf_{v_h \in CR_{h0}^1} | u -v_h |_{H^1(\mathbb{T}_h)} + \sup_{\varphi_h \in CR_{h0}^1} \frac{a_{0h}(u,\varphi_h) - (f,\varphi_h)}{| \varphi_h|_{H^1(\mathbb{T}_h)}}. \label{CR2}
\end{align}

The first term of the inequality \eqref{CR2} is estimated as follows. Using the CR interpolation error estimate \eqref{modify12}, we have, for any $u \in H^2(\Omega)$,
\begin{align}
\displaystyle
\inf_{v_h \in CR_{h0}^1} | u -v_h |_{H^1(\mathbb{T}_h)}
&\leq | u -I_h^{CR} u |_{H^1(\mathbb{T}_h)} \leq c H |u|_{H^2(\Omega)}. \label{CR3}
\end{align}

From the standard scaling argument, we have a consistency error inequality, e.g., see \cite[Lemma 3.36]{ErnGue04}.
\begin{lem}[Asymptotic Consistency]  \label{CR=lem6}
Let $u \in H^1_0(\Omega) \cap H^2(\Omega)$ be the solution of the homogeneous Dirichlet Poisson problem \eqref{intro1}. It then holds that 
\begin{align}
\displaystyle
\frac{a_{0h}(u,\varphi_h) - (f,\varphi_h)}{| \varphi_h|_{H^1(\mathbb{T}_h)}} \leq c \left( \sum_{T \in \mathbb{T}_h} \frac{h_T^4}{ (\min_{F \in \partial \mathbb{T}_h} \ell_F)^2} |u|^2_{H^2(T)} \right)^{1/2} \ \forall h, \ \forall \varphi_h \in  CR_{h0}^1, \label{CR4}
\end{align}
where $\partial \mathbb{T}_h$ denotes the set of all faces $F$ of $T \in \mathbb{T}_h$. Here, $\ell_F$ denotes the distance of the vertex of $T$ opposite to $F$ to the face. 
\end{lem}

\begin{pf*}
We follow \cite[Lemma 3.36]{ErnGue04}.

Let $\varphi_h \in  CR_{h0}^1$. Because $- \varDelta u = f$, we have
\begin{align*}
\displaystyle
a_{0h}(u,\varphi_h) - (f,\varphi_h)
&= \sum_{T \in \mathbb{T}_h} \int_T ( \nabla u \cdot \nabla \varphi_h - f \varphi_h ) dx \\
&= \sum_{T \in \mathbb{T}_h} \sum_{F \in \partial \mathbb{T}_h} \int_F (n_T \cdot \nabla ) u \varphi_h ds.
\end{align*}
Because each face $F$ of an element $T$ located inside $\Omega$ appears twice in the above sum, we have
\begin{align*}
\displaystyle
a_{0h}(u,\varphi_h) - (f,\varphi_h)
&= \sum_{T \in \mathbb{T}_h} \sum_{F \in \partial \mathbb{T}_h} \int_F (n_T \cdot \nabla ) u \left( \varphi_h - \overline{\varphi_h} \right) ds
\end{align*}
with the mean value
\begin{align*}
\displaystyle
\overline{\varphi_h} := \frac{1}{|F|} \int_F \varphi_h ds.
\end{align*}
Furthermore, we get
\begin{align*}
\displaystyle
a_{0h}(u,\varphi_h) - (f,\varphi_h)
&= \sum_{T \in \mathbb{T}_h} \sum_{F \in \partial \mathbb{T}_h} \int_F n_T \cdot \left(  \nabla u - \overline{\nabla u} \right)\left( \varphi_h - \overline{\varphi_h} \right) ds
\end{align*}
with the mean value
\begin{align*}
\displaystyle
n_T \cdot \overline{\nabla u} := \frac{1}{|F|} \int_F (n_T \cdot \nabla ) u ds.
\end{align*}
The Cauchy--Schwarz inequality yields
\begin{align*}
\displaystyle
a_{0h}(u,\varphi_h) - (f,\varphi_h)
&\leq  \sum_{T \in \mathbb{T}_h} \sum_{F \in \partial \mathbb{T}_h} \| \nabla u -\overline{\nabla u} \|_{L^2(F)^d} \| \varphi_h - \overline{\varphi_h} \|_{L^2(F)}.
\end{align*}
For $F \in \partial \mathbb{T}_h$, let $\widehat{T} \subset \mathbb{R}^d$ be the reference simplex and let $\Phi_T: \widehat{T} \to T$ be the corresponding affine transformation with Jacobian matrix $A_T$. Let $\widehat{F} = \Phi_T^{-1} ( F )$. Using the standard scaling argument and the trace theorem on the reference element, we have
\begin{align*}
\displaystyle
\| \varphi_h - \overline{\varphi_h} \|_{L^2(F)}
&\leq \left( \frac{|F|}{|\widehat{F}|} \right)^{1/2} \| \hat{\varphi}_h - \overline{\hat{\varphi}_h} \|_{L^2(\widehat{F})}
\leq c \left( \frac{|F|}{|\widehat{F}|} \right)^{1/2} \| \hat{\varphi}_h - \overline{\hat{\varphi}_h} \|_{H^1(\widehat{T})}.
\end{align*}
The Deny--Lions Lemma (see \cite[Lemma B.67]{ErnGue04}) implies
\begin{align*}
\displaystyle
 \| \hat{\varphi}_h - \overline{\hat{\varphi}_h} \|_{H^1(\widehat{T})}
 \leq c | \hat{\varphi}_h |_{H^1(\widehat{T})}.
\end{align*}
Using the standard scaling argument again, we obtain
\begin{align*}
\displaystyle
\| \varphi_h - \overline{\varphi_h} \|_{L^2(F)}
&\leq c \left( \frac{|F|}{|\widehat{F}|} \right)^{1/2}  | \hat{\varphi}_h |_{H^1(\widehat{T})} \\
&\leq  c \left( \frac{|F|}{|\widehat{F}|} \right)^{1/2} \| A_T \|_2 \left( \frac{|\widehat{T}|}{|{T}|} \right)^{1/2}  | {\varphi}_h |_{H^1({T})} \\
&\leq c \left( \frac{|F|}{|{T}|} \right)^{1/2} h_T | {\varphi}_h |_{H^1({T})} = c \left( \frac{d}{\ell_F} \right)^{1/2} h_T | {\varphi}_h |_{H^1({T})}.
\end{align*}
Here,  $\| A_T \|_2$ denotes the matrix $2$-norm as
\begin{align*}
\displaystyle
\| A_T \|_2 := \sup_{0 \neq x \in \mathbb{R}^d} \frac{|A_T x|}{|x|},
\end{align*}
where $| x | := ( \sum_{i=1}^d |x_i|^2)^{1/2}$ for $x \in \mathbb{R}^d$.

By analogous argument, we have
\begin{align*}
\displaystyle
 \| \nabla u -\overline{\nabla u} \|_{L^2(F)^d}
 &\leq c \left( \frac{d}{\ell_F} \right)^{1/2} h_T |u|_{H^2(T)}.
\end{align*}
We consequently get
\begin{align*}
\displaystyle
a_{0h}(u,\varphi_h) - (f,\varphi_h)
&\leq c \sum_{T \in \mathbb{T}_h} \sum_{F \in \partial \mathbb{T}_h} \frac{h_T^2}{\ell_F}  |u|_{H^2(T)} | {\varphi}_h |_{H^1({T})} \\
&\leq c \sum_{T \in \mathbb{T}_h} \frac{h_T^2}{\min_{F \in \partial \mathbb{T}_h} \ell_F}  |u|_{H^2(T)} | {\varphi}_h |_{H^1({T})} \\
&\leq c \left( \sum_{T \in \mathbb{T}_h} \frac{h_T^4}{(\min_{F \in \partial \mathbb{T}_h} \ell_F)^2}  |u|_{H^2(T)}^2  \sum_{T \in \mathbb{T}_h} | {\varphi}_h |_{H^1({T})}^2  \right)^{1/2},
\end{align*}
which leads to \eqref{CR4}.
\qed	
\end{pf*}

From \eqref{CR2}, \eqref{CR3} and \eqref{CR4}, we have
\begin{align*}
\displaystyle
| u - u_h^{CR} |_{H^1(\mathbb{T}_h)} 
\leq c H |u|_{H^2(\Omega)} + c \left( \sum_{T \in \mathbb{T}_h} \frac{h_T^4}{ (\min_{F \in \partial \mathbb{T}_h} \ell_F)^2} |u|^2_{H^2(T)} \right)^{1/2}.
\end{align*}

Since the order of the nonconforming term does not necessary becomes the order $H$, this inequality may be overestimated. 
\begin{ex*}
Let  $0 \< h_T \leq 1$. As examples, we consider two cases.
\begin{description}
  \item[(\Roman{lone})] When we use meshes including the tetrahedra $T$ with vertices $(0,0,0)^T$, $(h_T,0,0)^T$, $(0,h_T,0)^T$, and $(0,0,h_T^{\varepsilon})^T$, we have
 \begin{align*}
\displaystyle
 \frac{h_T^4}{(\min_{F \in \partial \mathbb{T}_h} \ell_F)^2} |u|^2_{H^2(T)}\leq c h_T^{2(2-\varepsilon)} |u|^2_{H^2(T)},
\end{align*}
where $1 \< \varepsilon \leq 2$. Since $H = \mathcal{O}(h)$, we get
\begin{align*}
\displaystyle
| u - u_h^{CR} |_{H^1(\mathbb{T}_h)} 
\leq c (h + h^{2 - \varepsilon}) |u|_{H^2(\Omega)}.
\end{align*}
  \item[(\Roman{ltwo})] When we use meshes including the tetrahedra $T$ with vertices $(0,0,0)^T$, $(h_T,0,0)^T$, $(0,h_T,0)^T$, and $(h_T^{\gamma},0,h_T^{\varepsilon})^T$, we have
 \begin{align*}
\displaystyle
 \frac{h_T^4}{(\min_{F \in \partial \mathbb{T}_h} \ell_F)^2} |u|^2_{H^2(T)} \leq c h_T^{2(2-\varepsilon)} |u|^2_{H^2(T)},
\end{align*}
where $1 \< \gamma \< \varepsilon \leq 1 + \gamma$ and  $1 \< \varepsilon \leq 2$. Since $H = \mathcal{O}(h^{1 + \gamma - \varepsilon})$, we get
\begin{align*}
\displaystyle
| u - u_h^{CR} |_{H^1(\mathbb{T}_h)} 
\leq c (h^{1 + \gamma - \varepsilon} + h^{2 - \varepsilon}) |u|_{H^2(\Omega)}.
\end{align*}
\end{description}
\end{ex*}

\subsection{Argument via the RT Interpolation Error}
To overcome the difficulty, we use the relation \eqref{pre21} in Lemma \ref{pre=lem2}, e.g., see also \cite{AcoDur99,LiuKik18}.

\begin{lem}[Asymptotic Consistency]  \label{CR=lem7}
We assume that $\Omega$ is convex. Let $\{ \mathbb{T}_h \}$ be a family of conformal meshes satisfying Assumption \ref{pre=ass1}. Let $u \in H^1_0(\Omega) \cap H^2(\Omega)$ be the solution of the homogeneous Dirichlet Poisson problem \eqref{intro1}. Then, there exists $c$, independent of $H$, such that
\begin{align}
\displaystyle
& \sup_{\varphi_h \in CR_{h0}^1} \frac{a_{0h}(u,\varphi_h) - (f,\varphi_h)}{| \varphi_h|_{H^1(\mathbb{T}_h)}} 
\leq c H \| f \|.  \label{CR5}
\end{align}
\end{lem}

\begin{pf*}
Using \eqref{pre21}, we have, for any $w_h \in RT_h^0$,
\begin{align*}
\displaystyle
\sup_{\varphi_h \in CR_{h0}^1} \frac{a_{0h}(u,\varphi_h) - (f,\varphi_h)}{| \varphi_h|_{H^1(\mathbb{T}_h)}} =
 \sup_{\varphi_h \in CR_{h0}^1} \frac{ (\nabla u - w_h , \nabla_h \varphi_h) - ( \div w_h + f, \varphi_h  ) }{|\varphi_h|_{H^1(\mathbb{T}_h)}}.
\end{align*}
We set $w_h := I_h^{RT} \nabla u$. From Lemma \ref{pre=lem1}, we get
\begin{align*}
\displaystyle
 \div ( I_h^{RT} \nabla u ) = \Pi_h^0 \div (\nabla u) = - \Pi_h^0 f.
\end{align*}
Furthermore, we have, for any $\varphi_h \in CR_{h0}^1$,
\begin{align*}
\displaystyle
  ( - \Pi_h^0 f + f,  \Pi_h^0 \varphi_h) &= 0.
\end{align*}
We thus obtain
\begin{align*}
\displaystyle
&(\nabla u - I_h^{RT} \nabla u, \nabla_h \varphi_h) - ( - \Pi_h^0 f + f, \varphi_h  ) \\
&= (\nabla u - I_h^{RT} \nabla u, \nabla_h \varphi_h) - ( - \Pi_h^0 f + f, \varphi_h  - \Pi_h^0 \varphi_h) \\
&\leq \| \nabla u - I_h^{RT} \nabla u \|_{L^2(\Omega)^d} |\varphi_h|_{H^1(\mathbb{T}_h)} + \| f -  \Pi_h^0 f \|  \| \varphi_h  - \Pi_h^0 \varphi_h \| \\
&\leq c H |u|_{H^2(\Omega)}  |\varphi_h|_{H^1(\mathbb{T}_h)} + c h \| f \| |\varphi_h|_{H^1(\mathbb{T}_h)}.
\end{align*}
\qed	
\end{pf*}

We consequently obtain the error estimate of the CR finite element method on anisotropic meshes.
\begin{thr} \label{CR=thr4}
We assume that $\Omega$ is convex. Let $\{ \mathbb{T}_h \}$ be a family of conformal meshes satisfying Assumption \ref{pre=ass1}. Let $u \in H^1_0(\Omega) \cap H^2(\Omega)$ be the solution of the homogeneous Dirichlet Poisson problem \eqref{intro1} with data $f \in L^2(\Omega)$. Let $u_h^{CR} \in CR_{h0}^1$ be the approximate solution of \eqref{CR1}. Then, there exists $c$, independent of $H$, such that
\begin{align}
\displaystyle
| u - u_h^{CR} |_{H^1(\mathbb{T}_h)} \leq c H \| f \|. \label{CR6}
\end{align}
\end{thr}

\begin{pf*}
Using \eqref{CR2}, \eqref{modify12} and \eqref{CR5}, we have
\begin{align*}
\displaystyle
| u - u_h^{CR} |_{H^1(\mathbb{T}_h)} 
&\leq 2 \inf_{v_h \in CR_{h0}^1} | u -v_h |_{H^1(\mathbb{T}_h)} + \sup_{\varphi_h \in CR_{h0}^1} \frac{a_{0h}(u,\varphi_h) - (f,\varphi_h)}{| \varphi_h|_{H^1(\mathbb{T}_h)}} \\
&\leq 2  | u -I_h^{CR} u |_{H^1(\mathbb{T}_h)} + c H \| f \|  \leq c H \| f \|,
\end{align*}
which leads to the estimate \eqref{CR6}.
\qed	
\end{pf*}

We next give the $L^2$ error estimate of the CR finite element method on anisotropic meshes, see also \cite{LiuKik07,LiuKik18,BreSco08}.

\begin{thr} \label{CR=thr5}
We assume that $\Omega$ is convex. Let $\{ \mathbb{T}_h \}$ be a family of conformal meshes satisfying Assumption \ref{pre=ass1}. Let $u \in H^1_0(\Omega) \cap H^2(\Omega)$ be the solution of the homogeneous Dirichlet Poisson problem \eqref{intro1} with data $f \in L^2(\Omega)$. Let $u_h^{CR} \in CR_{h0}^1$ be the approximate solution of \eqref{CR1}. Then, there exists $c$, independent of $H$, such that
\begin{align}
\displaystyle
\| u - u_h^{CR} \| \leq c H^2 \| f \|. \label{CR7}
\end{align}
\end{thr}

\begin{pf*}
We set $e_h :=  u - u_h^{CR} $. Let  $z \in H^2(\Omega) \cap H_0^1(\Omega)$ satisfy 
\begin{align}
\displaystyle
a_0(\varphi , z ) = ( \varphi , e_h) \quad \forall \varphi \in H_0^1(\Omega) \label{CR8}
\end{align}
and $z_h^{CR} \in CR_{h0}^1$ satisfy 
\begin{align}
\displaystyle
a_{0h}(\varphi_h , z_h^{CR}) = (\varphi_h , e_h) \quad \forall \varphi_h \in CR_{h0}^1. \label{CR9}
\end{align}
We then have
\begin{align}
\displaystyle
\| e_h \|^2 &= (e_h , e_h ) = a_{0h}(u , z )  - a_{0h}(u_h^{CR} , z_h^{CR}) \notag \\
&= a_{0h}(u -  u_h^{CR} , z - z_h^{CR} ) + a_{0h}( u -  u_h^{CR} ,  z_h^{CR}) + a_{0h}( u_h^{CR} ,  z - z_h^{CR} ) \notag \\
&= a_{0h}(u -  u_h^{CR} , z - z_h^{CR} ) \notag \\
&\quad + a_{0h}( u -  u_h^{CR} ,  z_h^{CR} - I_h^{CR} z) +   a_{0h}( u -  u_h^{CR} ,   I_h^{CR} z) \notag\\
&\quad + a_{0h}( u_h^{CR} - I_h^{CR} u ,  z - z_h^{CR} ) +a_{0h}( I_h^{CR} u ,  z - z_h^{CR} ).  \label{CR10}
\end{align}

Using Theorem \ref{CR=thr4}, the first term on the right hand side of \eqref{CR10} can be estimated as
\begin{align}
\displaystyle
 a_{0h}( u - u_h^{CR} ,  z- z_h^{CR})
 &\leq | u - u_h^{CR} |_{H^1(\mathbb{T}_h)}   | z- z_h^{CR} |_{H^1(\mathbb{T}_h)} \notag \\
 &\leq c H^2  \| f \| \|  e_h\|.  \label{CR11}
\end{align}
For the second and fourth terms on the right hand side of \eqref{CR10}, we have
\begin{align}
\displaystyle
& a_{0h}( u -  u_h^{CR} ,  z_h^{CR} - I_h^{CR} z) \notag \\
 &\quad =  a_{0h}( u -  u_h^{CR} ,  z_h^{CR} -  z) +  a_{0h}( u -  u_h^{CR} ,  z - I_h^{CR} z) \notag \\
 &\quad \leq  | u - u_h^{CR} |_{H^1(\mathbb{T}_h)} \left(  | z_h^{CR} - z |_{H^1(\mathbb{T}_h)} +  | z- I_h^{CR} z |_{H^1(\mathbb{T}_h)} \right) \notag \\
 &\quad \leq  c H^2  \| f \| \|  e_h\|, \label{CR12}
\end{align}
and, analogously,
\begin{align}
\displaystyle
a_{0h}( u_h^{CR} - I_h^{CR} u ,  z - z_h^{CR} ) 
&\leq c H^2 \| f \| \|  e_h\|. \label{CR13}
\end{align}
From \eqref{CR8}, \eqref{CR9} and \eqref{pre21}, we have
\begin{align*}
\displaystyle
&a_{0h}( u -  u_h^{CR} ,   I_h^{CR} z) \\
&\quad = a_{0h}( u ,   I_h^{CR} z) - a_{0h}( u_h^{CR} ,   I_h^{CR} z)  = (\nabla u , \nabla_h I_h^{CR} z) - (f ,  I_h^{CR} z) \\
&\quad = (\nabla u , \nabla_h I_h^{CR} z - \nabla z) - (f ,  I_h^{CR} z - z) + (\nabla u , \nabla z) - (f , z) \\
&\quad = (\nabla u - I_h^{RT} \nabla u , \nabla_h I_h^{CR} z - \nabla z) - (f + \div ( I_h^{RT} \nabla u ),  I_h^{CR} z - z).
\end{align*}
From Lemma \ref{pre=lem1} and $ \div ( I_h^{RT} \nabla u ) = - \Pi_h^0 f$, we have
\begin{align}
\displaystyle
&a_{0h}( u -  u_h^{CR} ,   I_h^{CR} z) \notag \\
&\quad =  (\nabla u - I_h^{RT} \nabla u , \nabla_h I_h^{CR} z - \nabla z) - (f - \Pi_h^0 f,  I_h^{CR} z - z) \notag\\
&\quad \leq \| \nabla u - I_h^{RT} \nabla u \|_{L^2(\Omega)^d} | I_h^{CR} z - z |_{H^1(\mathbb{T}_h)} + \| f - \Pi_h^0 f \| \|  I_h^{CR} z - z \|  \notag\\
&\quad \leq c H^2 \| f \| \|  e_h\|.  \label{CR14}
\end{align}

Analogously, from $ \div ( I_h^{RT} \nabla z ) = - \Pi_h^0 e_h$, we have
\begin{align}
\displaystyle
&a_{0h}( I_h^{CR} u ,  z - z_h^{CR} ) \notag \\
&\quad=  (\nabla_h I_h^{CR} u - \nabla u , \nabla z - I_h^{RT} \nabla z ) - (I_h^{CR} u - u , e_h + \div ( I_h^{RT} \nabla z)) \notag\\
&\quad \leq | I_h^{CR} u - u |_{H^1(\mathbb{T}_h)} \|  \nabla z - I_h^{RT} \nabla z \|_{L^2(\Omega)^d} + \|  I_h^{CR} u - u \| \| e_h  - \Pi_h^0 e_h \| \notag \\
&\quad \leq c H^2 \| f \| \|  e_h\|. \label{CR15}
\end{align}
Combining \eqref{CR10}, \eqref{CR11}, \eqref{CR12}, \eqref{CR13}, \eqref{CR14}, and \eqref{CR15}, we finally get 
\begin{align*}
\displaystyle
\| e_h \|^2 \leq c H^2 \| f \| \|  e_h\|, 
\end{align*}
which leads to the target estimate.
\qed	
\end{pf*}

\section{RT Finite Element Error Estimates}

\subsection{Dual mixed formulation of the Poisson problem}
The Poisson equation \eqref{intro1} $- \varDelta u =  - \div \nabla u = f$ can be written as the following system. Find $(\sigma,u):\Omega \to \mathbb{R}^{d} \times \mathbb{R}$ such that
\begin{subequations} \label{mix1}
\begin{align}
\displaystyle
\sigma - \nabla u &= 0  \quad \text{in $\Omega$}, \label{mix1a}\\
\div \sigma &=  - f  \quad \text{in $\Omega$}, \label{mix1b}\\
u &= 0 \quad \text{on $\partial \Omega$}. \label{mix1c}
\end{align}
\end{subequations}
We consider the following dual mixed formulation: Find $(\sigma , u) \in H(\div;\Omega) \times L^2(\Omega)$ such that
\begin{subequations} \label{mix2}
\begin{align}
\displaystyle
a(\sigma,v) + b(v,u) &= 0 \quad \forall v \in  H(\div;\Omega), \label{mix2a}\\
b(\sigma,q) &= - (f ,q) \quad \forall q \in L^2(\Omega), \label{mix2b}
\end{align}
\end{subequations}
where bilinear forms $a: H(\div;\Omega) \times H(\div;\Omega) \to \mathbb{R}$ and $b:H(\div;\Omega) \times L^2(\Omega) \to \mathbb{R}$ are defined by
\begin{align*}
\displaystyle
&a(\sigma,v) := (\sigma , v), \quad b(v,q) := (\div v , q).
\end{align*}
We set $X_0 := \{ v \in H(\div;\Omega)  ; \ b(v,q) = 0 \ \forall q \in L^2(\Omega) \}$. Because there exists a constant $c \> 0$ such that
\begin{align*}
\displaystyle
a(v,v) \geq c \| v \|^2_{H(\div;\Omega) } \quad \forall v \in X_0
\end{align*}
and the bilinear form $b(_\cdot , _ \cdot)$ satisfies the inf--sup condition
\begin{align}
\displaystyle
\inf_{0 \neq q \in L^2(\Omega)} \sup_{0 \neq v \in H(\div;\Omega)  } \frac{b(v,q)}{\| v \|_{H(\div;\Omega) } \| q \|}
\geq \beta_* \> 0, \label{mix3}
\end{align}
\eqref{mix2} is uniquely solvable; e.g., see \cite{GirRav86,BofBreFor13}. 

\subsection{RT Approximate Problem}
We consider the following RT approximate problem. Find $(\sigma_h^{RT} , u_h^{RT}) \in RT^0_{h} \times M_{h}^0$ such that
\begin{subequations} \label{mix4}
\begin{align}
\displaystyle
a(\sigma_h^{RT},v_h) + b(v_h,u_h^{RT}) &= 0, \quad \forall v_h \in  RT^0_{h}, \label{mix4a}\\
b(\sigma_h^{RT},q_h) &= - (f ,q_h), \quad \forall q_h \in M_{h}^0.\label{mix4b}
\end{align}
\end{subequations}
This setting is conforming because $RT^0_{h} \times M_{h}^0 \subset H(\div;\Omega) \times L^2(\Omega)$. It is given later that the discrete inf--sup condition
\begin{align*}
\displaystyle
\inf_{q_h \in M_{h}^0} \sup_{v_h \in RT^0_{h}} \frac{b(v_h , q_h)}{\| v_h \|_{H(\div;\Omega)} \| q_h \| } \geq c_* \> 0
\end{align*}
holds, where $c_*$ is a constant independent of $h$.

\subsection{Error Estimates of the RT Finite Element Approximation}
This section gives error estimates of the mixed finite element approximation \eqref{mix4}. We  emphasise that we do not impose the shape regularity condition and the maximum-angle condition for the mesh partition. That is, we assume that $\{ \mathbb{T}_h \}$ is a family of conformal meshes satisfying Assumption \ref{pre=ass1}.

\begin{lem} \label{mix=lem8}
Let $D \subset \mathbb{R}^d$ be a bounded domain. For any $g \in L^2(D)$, there exists $v \in H^1(D)^d$ such that
\begin{align}
\displaystyle
\div v = g \quad \text{in $D$} \label{mix5}
\end{align}
and
\begin{align}
\displaystyle
| v |_{H^1(D)^d} \leq  \| g \|_{L^2(D)}, \quad \| v \|_{L^2(\Omega)^d} \leq  C_P(D) \| g \|_{L^2(D)},  \label{mix6}
\end{align}
where $C_P(D)$ is the Poincar\'e constant.
\end{lem}

\begin{pf*}
The proof can be found in  \cite[Lemma 2.2]{BofBre08}.
\qed	
\end{pf*}

We next give the discrete inf--sup condition.
\begin{lem}[Discrete inf--sup condition] \label{mix=lem9}
If $ C_{G}^{RT} H \leq 1$, there exists a constant $c_*$, depending only on the Poincar\'e constant, such that   
\begin{align}
\displaystyle
\inf_{q_h \in M_{h}^0} \sup_{v_h \in RT^0_{h}} \frac{b(v_h , q_h)}{\| v_h \|_{H(\div;\Omega)} \| q_h \|} \geq c_* \> 0,\label{mix7}
\end{align}
where $ C_{G}^{RT}$ is the constant appearing in Corollary \ref{modify=coro3}.
\end{lem}

\begin{pf*}
Let $q_h \in  M_{h}^0$. From Lemma \ref{mix=lem8}, there exists $v \in H^1(\Omega)^d$ such that $\div v = q_h$ in $\Omega$, $| v |_{H^1(\Omega)^d} \leq   \| q_h \|$, and $\| v \|_{L^2(\Omega)^d} \leq C_P(\Omega) \| q_h \|$. 

By the Gauss theorem, we have
\begin{align*}
\displaystyle
\sum_{T \in \mathbb{T}_h} \int_{\partial T} v \cdot n_T ds = \sum_{T \in \mathbb{T}_h} \int_{T} \div v dx = \int_{\Omega} q_h dx.
\end{align*}
From the definition of the Raviart--Thomas interpolation, we conclude that
\begin{align*}
\displaystyle
 \int_{\Omega} \div ( I_T^{RT} v) p_h dx &= \sum_{T \in \mathbb{T}_h} p_h \int_T \div ( I_T^{RT} v) dx = \sum_{T \in \mathbb{T}_h} p_h \int_{\partial T} n_T \cdot (I_T^{RT} v) ds \\
&= \sum_{T \in \mathbb{T}_h} p_h \int_{\partial T} v \cdot n_T ds = \int_{\Omega} q_h p_h dx \quad \forall p_h \in M_h^0.
\end{align*}
Therefore,  it follows that $\div (I_h^{RT} v) = q_h$.

From the definitions, we have
\begin{align*}
\displaystyle
\| I_h^{RT} v \|^2_{H(\div;\Omega)}
&= \| I_h^{RT} v \|_{L^2(\Omega)^d}^2 + \| \div (I_h^{RT} v) \|^2 \\
&\leq 2 \| I_h^{RT} v - v \|_{L^2(\Omega)^d}^2 + 2 \| v \|_{L^2(\Omega)^d}^2 + \| q_h \|^2 \\
&\leq 2 ( C_{G}^{RT})^2 H^2 |v|^2_{H^1(\Omega)^d} + 2 C_P(\Omega)^2 \| q_h \|^2 + \| q_h \|^2 \\
&\leq \left( 3 + 2 C_P(\Omega)^2 \right) \| q_h \|^2. 
\end{align*}
We thus have
\begin{align*}
\displaystyle
 \sup_{v_h \in RT^0_{h}} \frac{b(v_h , q_h)}{\| v_h \|_{H(\div;\Omega)}}
& \geq \frac{b(I_h^{RT} v , q_h)}{\| I_h^{RT} v \|_{H(\div;\Omega)}}
 \geq \frac{1}{\left( 3 + 2 C_P(\Omega)^2 \right)^{1/2}} \frac{(q_h,q_h)}{\| q_h \|},
\end{align*}
and the proof of \eqref{mix7} is completed with $c_* := \left( 3 + 2 C_P(\Omega)^2 \right)^{- 1/2}$.
\qed
\end{pf*}

From the discrete equations \eqref{mix4} and their continuous counterpart \eqref{mix2}, we obtain the Galerkin orthogonality
\begin{subequations} \label{mix8}
\begin{align}
\displaystyle
a(\sigma - \sigma_h^{RT} , v_h) + b(v_h , u - u_h^{RT}) &= 0 \quad \forall v_h \in RT^0_{h}, \label{mix8a} \\
b(\sigma - \sigma_h^{RT} , q_h) &= 0 \quad \forall q_h \in M_{h}^0. \label{mix8b}
\end{align}
\end{subequations}
We then get the following C\'ea-lemma-type estimates with the help of \eqref{mix8} and the inf--sup condition \eqref{mix7}. 

\begin{thr} \label{mix=thr6}
Let $\sigma \in H^1(\Omega)^d$ and $\sigma_h^{RT} \in RT_h^0$ be the solutions of \eqref{mix1} and  \eqref{mix4}, respectively. We then have
\begin{align}
\displaystyle
\| \sigma - \sigma_h^{RT}\|_{L^2(\Omega)^d} \leq  \| \sigma - I_h^{RT} \sigma \|_{L^2(\Omega)^d}. \label{mix9}
\end{align}
Furthermore, let $(\sigma,u) \in H^1(\Omega)^d \times L^2(\Omega)$ and $(\sigma_h^{RT},u_h^{RT}) \in RT_h^0 \times M_h^0$ be the solutions of \eqref{mix1} and  \eqref{mix4}, respectively. Then, if $C_{G}^{RT} H \leq 1$,  it holds that
\begin{align}
\displaystyle
\| u - u_h^{RT} \|
\leq \| u - \Pi_h^0 u \| + c_*^{-1} \|  \sigma - \sigma_h^{RT}  \|_{L^2(\Omega)^d}. \label{mix11}
\end{align}
Here, $C_{G}^{RT}$ and $c_*$ are respectively the constants appearing in Corollary \ref{modify=coro3} and Lemma \ref{mix=lem9}.
\end{thr}

\begin{pf*}
The proof can be found in \cite[Lemma 3.7, Lemma 3.9]{BofBre08}.
\qed	
\end{pf*}

Using Theorem \ref{mix=thr6} and the interpolation error estimates of Corollary \ref{modify=coro1} and \ref{modify=coro3}, we thus have the error estimates of the mixed finite element approximation \eqref{mix4} on anisotropic meshes violating the maximum-angle condition.

\begin{thr} \label{mix=thr7}
let $(\sigma,u) \in H^1(\Omega)^d \times H^1(\Omega)$ and $(\sigma_h^{RT},u_h^{RT}) \in RT_h^0 \times M_h^0$ be the solutions of \eqref{mix1} and  \eqref{mix4}, respectively. Then, there exists a constant $c_1 \> 0$, independnt of  $\sigma$, $H$, and the geometric properties of $\mathbb{T}_h$,  such that
\begin{align}
\displaystyle
\| \sigma - \sigma_h^{RT}\|_{L^2(\Omega)^d}
&\leq c_1 H | \sigma |_{H^1(\Omega)^d}.  \label{mix10}
\end{align}
Furthermore, if $C_{G}^{RT} H \leq 1$, there exists a constant $c_2 \> 0$, depending on the discrete inf--sup condition but independent of $\sigma$, $u$, $h$, $H$, and the geometric properties of $\mathbb{T}_h$
\begin{align}
\displaystyle
\| u - u_h^{RT} \|
&\leq c_2 \left(h |u|_{H^1(\Omega)} + H |\sigma|_{H^1(\Omega)^d} \right). \label{mix12}
\end{align}
Here, $C_{G}^{RT}$ is the constant appearing in Corollary \ref{modify=coro3}.
\end{thr}

\section{Relationship between the RT and CR Finite Element Approximation}
This section shows the relationship between the RT and CR problems. Find $(\bar{\sigma}_h^{RT} , \bar{u}_h^{RT}) \in RT^0_{h} \times M_{h}^0$ such that
\begin{subequations} \label{RTCR1}
\begin{align}
\displaystyle
a(\bar{\sigma}_h^{RT},v_h) + b(v_h,\bar{u}_h^{RT}) &= 0 \quad \forall v_h \in  RT^0_{h}, \label{RTCR1a}\\
b(\bar{\sigma}_h^{RT},q_h) &= - (\Pi_h^0 f ,q_h) \quad \forall q_h \in M_{h}^0 \label{RTCR1b}
\end{align}
\end{subequations}
and find $\bar{u}_h^{CR} \in CR_{h0}^1$ such that
\begin{align}
\displaystyle
a_{0h}(\bar{u}_h^{CR} , \varphi_h) = (\Pi_h^0 f , \varphi_h) \quad \forall \varphi_h \in CR_{h0}^1. \label{RTCR2}
\end{align}
Here, \eqref{RTCR2} is the CR approximation of the Poisson equation
\begin{align}
\displaystyle
- \varDelta \bar{u} = \Pi_h^0 f \quad \text{in $\Omega$}, \quad \bar{u} = 0 \quad \text{on $\partial \Omega$}. \label{RTCR3}
\end{align}

In the case of $d=2$, it is well known that there exists a relationship between $(\bar{\sigma}_h^{RT} , \bar{u}_h^{RT}) $ and $\bar{u}_h^{CR}$ introduced by Marini; for example, \cite{Mar85}. See also \cite{LiuKik07,KikSai16,LiuKik18}. We here show the relation in the three dimensional case. 

\begin{figure}
\includegraphics[width=7cm]{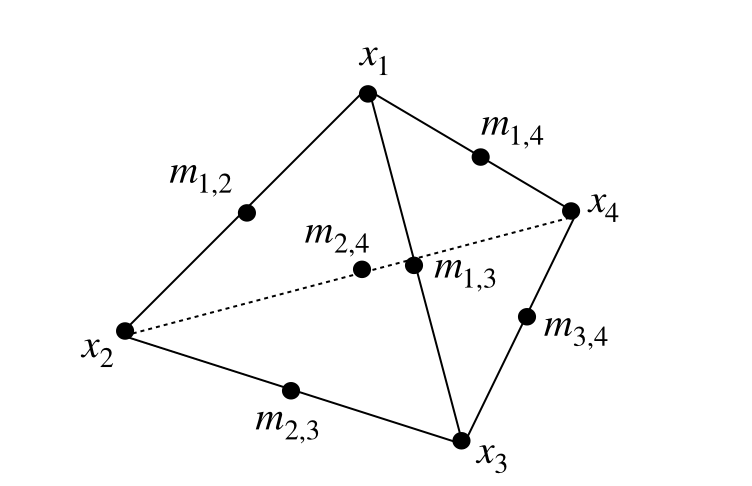}
\caption{Tetrahedron}
\label{fig1}	
\end{figure}

Let us consider a tetrahedron $T \subset \mathbb{R}^3$ such as that in Figure \ref{fig1}. Let $x_i$ ($i=1,2,3,4$) be the vertices and $m_{i,j}$ the midpoints of edges of the tetrahedron; that is, $m_{i,j} := \frac{1}{2} (x_i + x_j)$. Furthermore, for $1 \leq i \leq 4$, let $F_i$ be the face of the tetrahedron opposite $x_i$. Then, by simple calculation, we find the equality
\begin{align*}
\displaystyle
L := \sum_{i=1}^4 |x_i - x_T|^2 =  |m_{1,4} - m_{2,3}|^2 +  |m_{1,3} - m_{2,4}|^2 +  |m_{1,2} - m_{3,4}|^2,
\end{align*}
holds, where $x_T$ is the barycentre of $T$ such that $x_T := \frac{1}{4} \sum_{i=1}^4 x_i$. 

We present a quadrature scheme over a simplex $T \subset \mathbb{R}^3$ (e.g., \cite[p.307]{Str71}) that is easily conformed. 
\begin{lem} \label{RTCR=lem10}
For any $f \in \mathcal{C}^0(T)$, the quadrature scheme
\begin{align*}
\displaystyle
\int_T f(x) dx &\sim - \frac{|T|}{20} \sum_{i=1}^4 f(x_i) + \frac{|T|}{5} \sum_{1 \leq i \textless j \leq 4} f(m_{i,j})
\end{align*}
is exact for polynomials of degree less than or equal to $2$;
\begin{align}
\displaystyle
\int_T f(x) dx +  \frac{|T|}{20} \sum_{i=1}^4 f(x_i) - \frac{|T|}{5} \sum_{1 \leq i \< j \leq 4}f(m_{i,j}) = 0 \quad \forall f \in \mathcal{P}^2(T). \label{RTCR4}
\end{align}
\end{lem}

Define the function $\varphi_T$ by
\begin{align}
\displaystyle
\varphi_T(x) &:= 
\begin{cases} \label{RTCR5}
L - 12 |x - x_T|^2, \quad \text{on $T$},\\
0, \quad \text{otherwise}.
\end{cases}
\end{align}
We then have the following lemma.

\begin{lem}  \label{RTCR=lem11}
It holds that
\begin{align}
\displaystyle
\frac{1}{|F_i|} \int_{F_i} \varphi_T (x) ds &= 0, \quad i=1,2,3,4, \label{RTCR6} \\
\frac{1}{|T|} \int_T \varphi_T (x) dx &= \frac{2}{5} L, \label{RTCR7} \\
\frac{1}{|T|} \int_T | \nabla \varphi_T (x)|^2 dx &=  \frac{144}{5} L.\label{RTCR8}
\end{align}
\end{lem}

\begin{pf*}
From second-order three-point numerical integration over $F_1$, 
\begin{align*}
\displaystyle
\int_{F_1} f(x) ds &= \frac{|F_1|}{3} \left( f (m_{2,3}) + f (m_{3,4}) + f (m_{2,4}) \right)  \quad \forall f \in \mathcal{P}^2(T),
\end{align*}
we have
\begin{align*}
\displaystyle
&\frac{1}{|F_1|} \int_{F_1} \varphi_T (x) ds \\
&= \frac{1}{3} \left( \varphi_T (m_{2,3}) + \varphi_T (m_{3,4}) + \varphi_T (m_{2,4}) \right) \\
&= \frac{1}{3} \left( 3 L - 12 \left(  |m_{2,3} - x_T|^2 + |m_{3,4} - x_T|^2 + |m_{2,4} - x_T|^2 \right) \right) \\
&= \frac{1}{3} \left( 3 L - \frac{12}{4} \left(  |m_{2,3} - m_{1,4}|^2 + |m_{3,4} - m_{1,2}|^2 + |m_{2,4} - m_{1,3}|^2 \right) \right) = 0,
\end{align*}
which leads to \eqref{RTCR6}.

Next, using \eqref{RTCR4}, we have
\begin{align*}
\displaystyle
&\frac{1}{|T|} \int_T \varphi_T (x) dx \\
&= - \frac{1}{20} \sum_{i=1}^4 \varphi_T(x_i) + \frac{1}{5} \sum_{1 \leq i \< j \leq 4}\varphi_T(m_{i,j}) \\
&= - \frac{1}{20} \left( 4 L - 12 \sum_{i=1}^4 |x_i - x_T|^2 \right) + \frac{1}{5} \left( 6L - 12 \sum_{1 \leq i \< j \leq 4} |m_{i,j} - x_T|^2 \right) \\
&= \frac{2}{5} L,
\end{align*}
which leads to \eqref{RTCR7}. We here used
\begin{align*}
\displaystyle
&\sum_{1 \leq i \< j \leq 4} |m_{i,j} - x_T|^2  \\
&=  |m_{1,2} - x_T|^2 +  |m_{1,3} - x_T|^2 +  |m_{1,4} - x_T|^2 \\
&\quad +  |m_{2,3} - x_T|^2 +  |m_{2,4} - x_T|^2 +  |m_{3,4} - x_T|^2 \\
&= \frac{1}{4} \left( 2 |m_{1,2} - m_{3,4}|^2 + 2 |m_{1,3} - m_{2,4}|^2 +2  |m_{1,4} - m_{2,3}|^2  \right) = \frac{L}{2}.
\end{align*}
We similarly obtain
\begin{align*}
\displaystyle
&\frac{1}{|T|} \int_T | \nabla \varphi_T (x)|^2 dx \\
&= \frac{24^2}{|T|} \int_T |x - x_T|^2 dx \\
&= - \frac{24^2}{20} \sum_{i=1}^4 |x_i - x_T|^2 + \frac{24^2}{5} \sum_{1 \leq i \< j \leq 4} |m_{i,j} - x_T|^2 = \frac{144}{5} L,
\end{align*}
which leads to \eqref{RTCR8}.
\qed
\end{pf*}

We set the bubble space $B_h$ by
\begin{align}
\displaystyle
B_h := \{ b_h \in L^2(\Omega); \ b_h |_T \in \Span \{ \varphi_T\}, \ \forall T \in \mathbb{T}_h \}. \label{RTCR9}
\end{align}
Then, for any $\psi_h \in CR_{h0}^1$ and $b_h \in B_h$, because one writes $b_h|_T = c_b\varphi_T$ for $c_b \in \mathbb{R}$, it holds that
\begin{align*}
\displaystyle
(\nabla_h \psi_h, \nabla_h b_h)  &= \sum_{T \in \mathbb{T}_h} c_b \int_T \nabla \psi_h \cdot \nabla \varphi_T dx \\
&= \sum_{T \in \mathbb{T}_h} c_b \left\{ \sum_{F \subset \partial T} (n_F \cdot \nabla \psi_h ) \int_F \varphi_T ds - \int_T \varDelta \psi_h \varphi_T dx \right\} = 0.
\end{align*}
We here used the facts that \eqref{RTCR6}, $n_F \cdot \nabla \psi_h$ is constant on $F$, and $\varDelta \psi_h = 0$ on $T$. That is to say, two finite element spaces $CR_{h0}^1$ and $B_h$ are orthogonal to  each other. 

Furthermore, we define the finite element space $X_h^{bCR}$ by
\begin{align}
\displaystyle
X_h^{bCR} := CR_{h0}^1 + B_h = \{ \psi_h + b_h ; \ \psi_h \in CR_{h0}^1, \ b_h \in B_h  \}. \label{RTCR10}
\end{align}
We consider the following finite element problem. Find $u_h^{bCR} \in X_h^{bCR}$ such that
\begin{align}
\displaystyle
a_{0h}(u_h^{bCR} , \varphi_h) = (\nabla_h u_h^{bCR} , \nabla_h \varphi_h) = (\Pi_h^0 f , \varphi_h) \quad \forall \varphi_h \in X_{h}^{bCR}.\label{RTCR11}
\end{align}
The solution $u_h^{bCR} \in X_h^{bCR}$ is then decomposed as $u_h^{bCR} = \bar{u}_h^{CR} + b_h$ with $\bar{u}_h^{CR} \in  CR_{h0}^1$ and $b_h\in  B_h$. Note that $\bar{u}_h^{CR}$ and $b_h$  respectively satisfy \eqref{RTCR2} and the equation
\begin{align}
\displaystyle
a_{0h}(b_h , c_h) = (\nabla_h b_h , \nabla_h c_h) = (\Pi_h^0 f , c_h) \quad \forall c_h \in  B_h. \label{RTCR12}
\end{align}
On each element $T \in \mathbb{T}_h$, \eqref{RTCR12} has the form
\begin{align*}
\displaystyle
\gamma_T \int_T \nabla \varphi_T \cdot \nabla \varphi_T dx = \int_T \Pi_T^0 f \varphi_T dx, \quad \gamma_T \in \mathbb{R}. 
\end{align*}
From \eqref{RTCR7} and \eqref{RTCR8}, we have
\begin{align}
\displaystyle
\gamma_T = \frac{1}{72} \Pi_T^0 f \quad \forall T \in \mathbb{T}_h. \label{RTCR13}
\end{align}
\begin{thr}  \label{RTCR=thr8}
Let ${u}_h^{bCR} \in X_{h}^{bCR}$ be the solution of \eqref{RTCR11} and $(\bar{\sigma}_h^{RT},\bar{u}_h^{RT}) \in RT^0_{h} \times M_{h}^0$ the solution of \eqref{RTCR1}. We then have $\nabla_h u_h^{bCR} \in RT^0_{h}$ and
\begin{align}
\displaystyle
\bar{\sigma}_h^{RT} &= \nabla u_h^{bCR} \quad\forall T \in \mathbb{T}_h,  \label{RTCR14}\\
\bar{u}_h^{RT} &= \Pi_T^0  u_h^{bCR}  \quad\forall T \in \mathbb{T}_h.\label{RTCR15}
\end{align}
\end{thr}

\begin{pf*}
The proof can be found in \cite{HuMa15}.
\qed	
\end{pf*}

From Theorem \ref{RTCR=thr8}, for $d=3$, the following lemma holds.
\begin{lem}  \label{RTCR=lem12}
Let $\overline{u}_h^{CR} \in CR_{h0}^1$ be the solution of \eqref{RTCR2} and $(\bar{\sigma}_h^{RT},\bar{u}_h^{RT}) \in RT^0_{h} \times M_{h}^0$ be the solution of \eqref{RTCR1}. We then have the relationships
\begin{align}
\displaystyle
\bar{\sigma}_h^{RT}|_T &= \nabla \bar{u}_h^{CR} - \frac{1}{3} \Pi_T^0 f (x - x_T) \quad \forall T \in \mathbb{T}_h, \label{RTCR16} \\
\bar{u}_h^{RT}|_T &= \Pi_T^0  \bar{u}_h^{CR} + \frac{1}{180} \Pi_T^0 f  \sum_{i=1}^4 |x_i - x_T|^2  \quad \forall T \in \mathbb{T}_h. \label{RTCR17}
\end{align}
\end{lem}

Using relationship between the RT and CR finite element methods, we have the error estimate of the CR finite element approximation with the bubble function.

\begin{lem} \label{RTCR=lem13}
We assume that $\Omega$ is convex. Let $\{ \mathbb{T}_h \}$ be a family of conformal meshes satisfying Assumption \ref{pre=ass1}. Let $\bar{u} \in  H_0^1(\Omega) \cap H^2(\Omega)$ be the solution of \eqref{RTCR3} and $u_h^{bCR} \in X_{h}^{bCR}$ be the solution of the CR problem \eqref{RTCR11}. There then exists  a constant $c \> 0$ independent of $\bar{u}$, $h$, $H$ and the geometric properties of $\mathbb{T}_h$ such that
\begin{align}
\displaystyle
| \bar{u} - u_h^{bCR} |_{H^1(\mathbb{T}_h)} \leq  c H \| \Pi_h^0 f \|. \label{RTCR18}
\end{align}
\end{lem}

\begin{pf*}
Let $(\bar{\sigma}_h^{RT},\bar{u}_h^{RT}) \in RT^0_{h} \times M_{h}^0$ be the solution of \eqref{RTCR1}. From Theorem \ref{RTCR=thr8}, it holds that $\nabla_h u_h^{bCR} \in RT^0_{h}$ and $\bar{\sigma}_h^{RT} = \nabla_h u_h^{bCR}$. Setting $\bar{\sigma} := \nabla \bar{u} \in H^1(\Omega)^d$, we then have, using inequality \eqref{mix10}, that
\begin{align*}
\displaystyle
| \bar{u} - u_h^{bCR} |_{H^1(\mathbb{T}_h)}
&=  \left( \sum_{T \in \mathbb{T}_h} \| \bar{\sigma} - \bar{\sigma}_h^{RT} \|^2_{L^2(T)^d} \right)^{1/2} \\
&\leq c H |\bar{\sigma}|_{H^1(\Omega)^d} = c H |\bar{u}|_{H^2(\Omega)} \leq c H \| \Pi_h^0 f \|.
\end{align*}
\qed
\end{pf*}

\section{Numerical Results}
This section presents results of numerical examples. Let $\Omega := (0,1)^3$. Let $u_h^L$ and $u_h^{CR}$ be the $\mathcal{P}^1$-Lagrange and $\mathcal{P}^1$-CR finite element solutions, respectively, for the model problem
\begin{align*}
\displaystyle
- \varDelta u &= 2y(1-y)z(1-z) + 2x(1-x)z(1-z) + 2x(1-x)y(1-y) \quad \text{in $\Omega$}, \\
u &= 0 \quad \text{on $\partial \Omega $},
\end{align*}
which is the exact solution $u = x(1-x)y(1-y)z(1-z)$. 

\begin{figure}[htbp]
  \includegraphics[width=7cm]{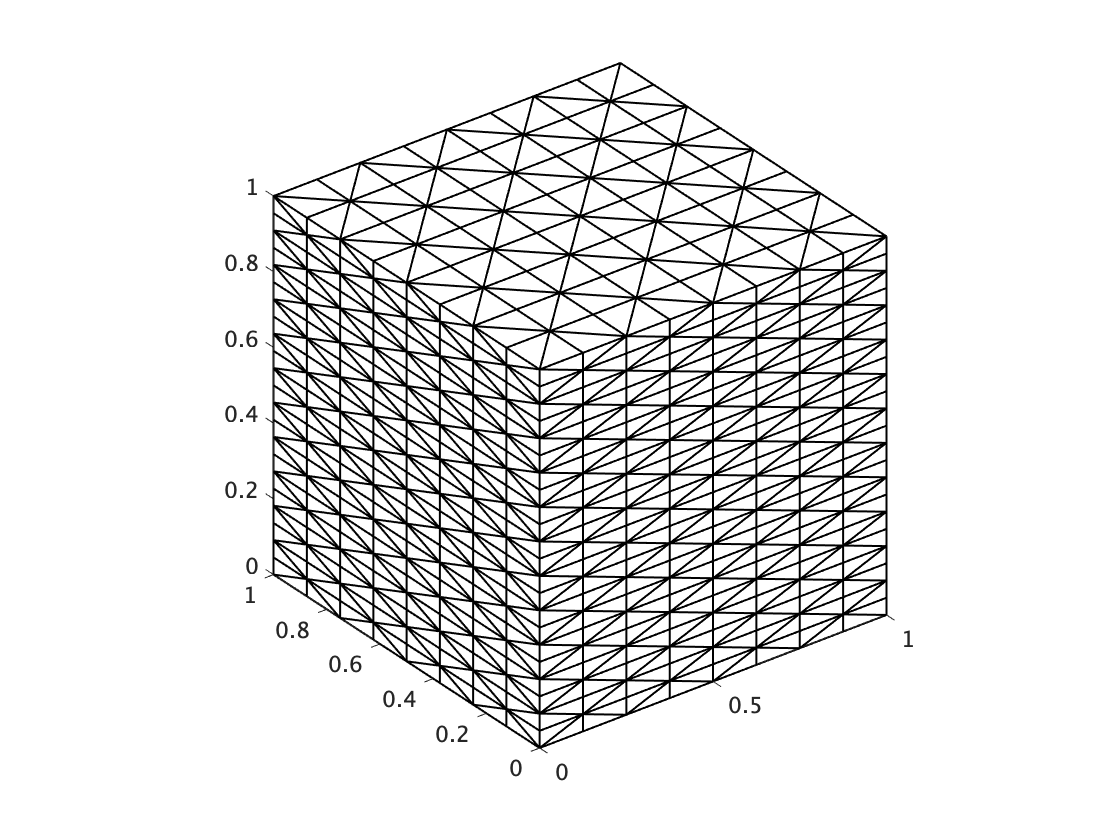}
  \caption{Mesh: $M=8$, $N=22$}
   \label{fig2}
   \vspace{0.5cm}
   \includegraphics[width=8cm]{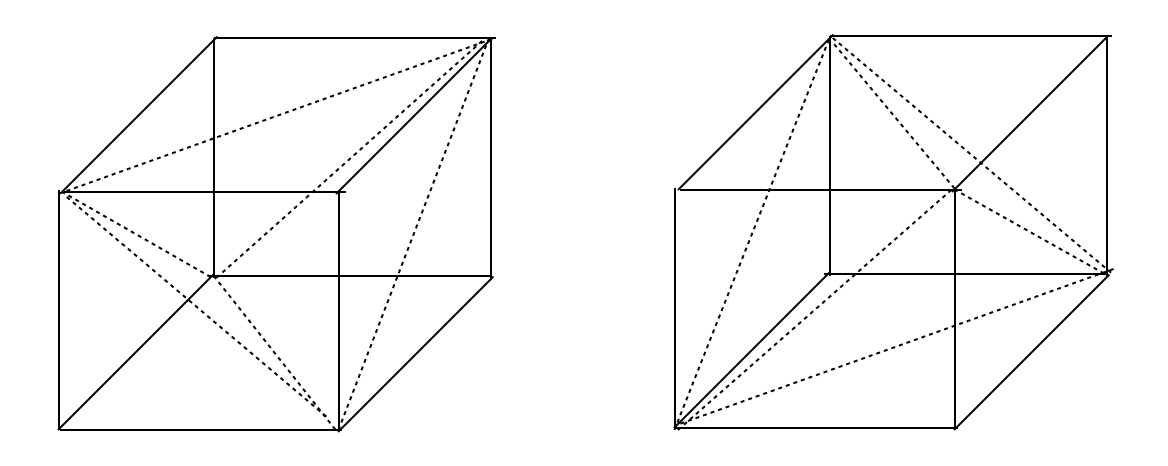}
  \caption{Elements}
   \label{fig3}
\end{figure}

Let $M$ be the division number of each side of the bottom face and $N$ the division number of the height of $\Omega$ with $N \sim M^{\gamma}$ (see Fig. \ref{fig2}). There are two elements as shown in Fig. \ref{fig3}.

If an exact solution $u$ is known, the error $e_h := u - u_h$ and $e_{h/2} := u - u_{h/2}$ are computed numerically for two mesh sizes $h$ and $h/2$. The convergence indicator $r$ is defined by
\begin{align*}
\displaystyle
r = \frac{1}{\log(2)} \log \left( \frac{\| e_h \|_X}{\| e_{h/2} \|_X} \right).
\end{align*}

We set $h := \frac{1}{M}$. The parameter $H$ is then $H = \mathcal{O}( h^{2 - \gamma} )$. We compute the convergence order with respect to $H_0^1$ and $L^2$ norms defined by
\begin{align*}
\displaystyle
&Err_h^L(H^1) := \frac{| u - u_h^L |_{H^1(\Omega)}}{\| \varDelta u \|}, \quad  Err_h^{L}(L^2) := \frac{\| u - u_h^L \|}{\| \varDelta u \|}, \\
&Err_h^{CR}(H^1) := \frac{| u - u_h^{CR} |_{H^1(\mathbb{T}_h)}}{\| \varDelta u \|}, \quad  Err_h^{CR}(L^2) := \frac{\| u - u_h^{CR} \|}{\| \varDelta u \|},
\end{align*}
for three cases: $\gamma = 1.5$, $\gamma = 1.9$ and $\gamma = 2.0$. In order to compute the above norms, we use the five-order fifteen-point numerical integration introduced in \cite{Kea86}. The results are give in Table \ref{table1}, Table \ref{table2} when $\gamma = 1.5$, Table \ref{table3}, Table \ref{table4} when $\gamma = 1.9$, and Table \ref{table5}, Table \ref{table6} when $\gamma = 2.0$. Further, $N_p^L$ and  $N_p^{CR}$ denote respectively the degrees of freedom for the $\mathcal{P}^1$-Lagrange finite element and the $\mathcal{P}^1$-CR finite element.

\begin{table}[htb]
\caption{Error of the $\mathcal{P}^1$-Lagrange finite element solution ($\gamma = 1.5$)}
\centering
\begin{tabular}{l|l|l|l|l|l|l|l|l} \hline
$M$& $N$ &  $h$ & $H$ & $N_p^L$ & $Err_h^L(H^1)$ & $r$& $Err_h^L(L^2)$  & $r$  \\ \hline \hline
4& 8 & 2.50e-01 & 5.00e-01  &225 &  1.2043e-01&   & 9.5321e-03&  \\
8& 22 & 1.25e-01 & 3.54e-01 & 1,863  & 7.0318e-02   & 0.78   & 3.1646e-03    & 1.59   \\
16& 64 & 6.25e-02 & 2.50e-01 & 18,785 & 4.4662e-02  &  0.65 & 1.2570e-03   & 1.33 \\
32& 182 &  3.13e-02 & 1.77e-01 & 199,287 & 2.9479e-02   & 0.60  & 5.4477e-04   & 1.21 \\
\hline
\end{tabular}
\label{table1}
\end{table}

\begin{table}[htb]
\caption{Error of the $\mathcal{P}^1$-CR finite element solution ($\gamma = 1.5$)}
\centering
\begin{tabular}{l|l|l|l|l|l|l|l|l} \hline
$M$& $N$ &  $h $ & $H$ & $N_p^{CR}$ & $Err_h^{CR}(H^1)$ & $r$ & $Err_h^{CR}(L^2)$  & $r$  \\ \hline \hline
4& 8 &2.50e-01  &5.00e-01 & 1,440 &8.2569e-02   &   & 3.8242e-03    &  \\
8&  22 & 1.25e-01 & 3.54e-01 & 14,912 & 4.0629e-02  & 1.02   & 8.8356e-04    & 2.11 \\
16& 64 & 6.25e-02 & 2.50e-01 & 168,448 &  2.0042e-02 & 1.02  & 2.0485e-04    & 2.11 \\
32& 182 &   3.13e-02  &1.77e-01 & 1,889,024  &  9.9579e-03 & 1.01  & 4.8960e-05   & 2.07 \\
\hline
\end{tabular}
\label{table2}
\end{table}

\begin{table}[htb]
\caption{Error of the $\mathcal{P}^1$-Lagrange finite element solution ($\gamma = 1.9$)}
\centering
\begin{tabular}{l|l|l|l|l|l|l|l|l} \hline
$M$& $N$ &  $h $ & $H$ & $N_p^L$ & $Err_h^L(H^1)$ & $r$ & $Err_h^L(L^2)$  & $r$  \\ \hline \hline
4& 14  & 2.50e-01 & 8.71e-01  & 345 &  1.4873e-01 &   & 1.4032e-02   &  \\
8& 52 & 1.25e-01 & 8.12e-01 & 4,293 & 1.2167e-01  &  0.29 &  9.3061e-03  & 0.59 \\
16& 194 & 6.25e-02 &  7.58e-01 & 56,355 & 1.0919e-01  & 0.16  &7.4989e-03    & 0.31\\
32& 724 & 3.13e-02 & 7.07e-01 & 789,525  & 1.0128e-01  &  0.11 & 6.4558e-03    & 0.22 \\
\hline
\end{tabular}
\label{table3}
\end{table}

\begin{table}[htb]
\caption{Error of the $\mathcal{P}^1$-CR finite element solution ($\gamma = 1.9$)}
\centering
\begin{tabular}{l|l|l|l|l|l|l|l|l} \hline
$M$& $N$ &  $h $ & $H$ & $N_p^{CR}$ & $Err_h^{CR}(H^1)$ & $r$ & $Err_h^{CR}(L^2)$  & $r$ \\ \hline \hline
4& 14 & 2.50e-01 & 8.71e-01 & 2,496 & 7.9756e-02 &   & 3.2993e-03   &  \\
8& 52 & 1.25e-01 & 8.12e-01 & 35,072  & 3.9708e-02   & 1.01   & 7.7177e-04  & 2.10 \\
16& 194 & 6.25e-02 & 7.58e-01 & 509,568  & 1.9814e-02  & 1.00  &   1.8781e-04 &2.04 \\
32& 724 &  3.13e-02 & 7.07e-01 & 7,508,480 &  9.9003e-03  &  1.00 & 4.6546e-05    & 2.01 \\
\hline
\end{tabular}
\label{table4}
\end{table}

\begin{table}[htb]
\caption{Error of the $\mathcal{P}^1$-Lagrange finite element solution ($\gamma = 2.0$)}
\centering
\begin{tabular}{l|l|l|l|l|l|l|l|l} \hline
$M$& $N$ &  $h$ & $H$ & $N_p^L$ & $Err_h^L(H^1)$ & $r$ & $Err_h^L(L^2)$  & $r$ \\ \hline \hline
4& 16  & 2.50e-01 & 1.00   & 425 & 1.5862e-01  &   & 1.5909e-02   &  \\
8& 64 & 1.25e-01 & 1.00  & 5,265 & 1.4079e-01  &  0.17 & 1.2472e-02   & 0.35 \\
16& 256 & 6.25e-02  & 1.00 & 74,273 & 1.3597e-01  & 0.05  & 1.1646e-02    & 0.10 \\
32& 1,024 & 3.13e-02 & 1.00 & 1,116,225  & 1.3474e-01   & 0.01  & 1.1442e-02   & 0.03 \\
\hline
\end{tabular}
\label{table5}
\end{table}

\begin{table}[htb]
\caption{Error of the $\mathcal{P}^1$-CR finite element solution ($\gamma = 2.0$)}
\centering
\begin{tabular}{l|l|l|l|l|l|l|l|l} \hline
$M$& $N$ &  $h $ & $H$& $N_p^{CR}$ & $Err_h^{CR}(H^1)$ & $r$ & $Err_h^{CR}(L^2)$  & $r$  \\ \hline \hline
4& 16  & 2.50e-01 & 1.00 & 2,848 & 7.9473e-02  &   &   3.2264e-03 &  \\
8& 64 & 1.25e-01 & 1.00  & 43,136  & 3.9647e-02  & 1.00  &  7.6153e-04  & 2.08  \\
16& 256  & 6.25e-02  & 1.00  & 672,256 & 1.9803e-02   & 1.00  & 1.8680e-04    & 2.03  \\
32& 1,024  & 3.13e-02 & 1.00 & 10,618,880  &  9.8984e-03 & 1.00  &  4.6458e-05  & 2.01 \\
\hline
\end{tabular}
\label{table6}
\end{table}

Observing the numerical results, the convergence indicators $r$ in each norms are respectively
\begin{align*}
\displaystyle
&| u - u_h^L |_{H^1(\Omega)} = \mathcal{O}(H), \quad \| u - u_h^L \| = \mathcal{O}(H^2), \\
&| u - u_h^{CR} |_{H^1(\mathbb{T}_h)} = \mathcal{O}(h), \quad \| u - u_h^{CR} \| = \mathcal{O}(h^2),
\end{align*}
where $H = \mathcal{O}( h^{2 - \gamma} )$. Meanwhile, the theoretical results are as follows:
\begin{align*}
\displaystyle
&| u - u_h^L |_{H^1(\Omega)} = \mathcal{O}(H), \quad \| u - u_h^L \| = \mathcal{O}(H^2), \\
&| u - u_h^{CR} |_{H^1(\mathbb{T}_h)} = \mathcal{O}(H), \quad \| u - u_h^{CR} \| = \mathcal{O}(H^2),
\end{align*}
if $\Omega$ is convex and $u \in H^2(\Omega) \cap H_0^1(\Omega)$. In this numerical examples, the CR finite element approximation is superior to the Lagrange finite element approximation on this anisotropic meshes. The theoretical explanation of this point is still open.


%
%

\begin{acknowledgements}
This work was supported by JSPS KAKENHI Grant Number \\ JP16H03950. We would like to thank the anonymous referee for the valuable comments.
\end{acknowledgements}

%
%



\end{document}